\documentclass[11pt,twoside,a4paper]{article}
\usepackage{amsmath,amssymb,amsthm}

\usepackage{graphicx,psfrag}

\newtheorem{thm}{Theorem}[section]
\newtheorem{lem}[thm]{Lemma}
\newtheorem{pro}[thm]{Proposition}
\newtheorem{cor}[thm]{Corollary}

\theoremstyle{definition}
\newtheorem{Def}[thm]{Definition}
\newtheorem{exa}[thm]{Example}

\theoremstyle{remark}
\newtheorem{rem}[thm]{Remark}

\newcommand{\R}{\mathbb{R}}
\newcommand{\Z}{\mathbb{Z}}
\newcommand{\N}{\mathbb{N}}
\renewcommand{\H}{\mathbb{H}}

\newcommand{\cL}{\mathcal{L}}

\newcommand{\cU}{\mathcal{U}}

\newcommand{\al}{\alpha}
\newcommand{\be}{\beta}
\newcommand{\ga}{\gamma}
\newcommand{\Ga}{\Gamma}
\newcommand{\de}{\delta}

\newcommand{\om}{\omega}
\newcommand{\Om}{\Omega}
\newcommand{\si}{\sigma}

\newcommand{\ka}{\kappa}
\newcommand{\la}{\lambda}
\newcommand{\La}{\Lambda}
\renewcommand{\phi}{\varphi}

\newcommand{\dist}{\operatorname{dist}}
\newcommand{\diam}{\operatorname{diam}}

\newcommand{\asdim}{\operatorname{asdim}}
\newcommand{\hypdim}{\operatorname{hypdim}}

\newcommand{\mesh}{\operatorname{mesh}}
\newcommand{\cdim}{\operatorname{cdim}}

\newcommand{\lv}{\operatorname{lv}}

\newcommand{\es}{\emptyset}
\renewcommand{\d}{\partial}
\newcommand{\di}{\d_{\infty}}
\newcommand{\set}[2]{\{#1:\,\text{#2}\}}
\newcommand{\sm}{\setminus}
\newcommand{\sub}{\subset}

\newcommand{\un}{\underline}
\newcommand{\ov}{\overline}

\newcommand{\wh}{\widehat}

\begin{document}

\title{A product of trees as universal space for hyperbolic groups}
\author{Sergei Buyalo\footnote{Supported by RFFI Grant
05-01-00939, Grant NSH-1914.2003.1 and SNF Grant 20-668 33.01}
\ \& Viktor Schroeder\footnote{Supported by Swiss National Foundation}}

\date{}
\maketitle

\begin{abstract}
We show that every Gromov hyperbolic group
$\Ga$
admits a quasi-isometric embedding into the product of
$(n+1)$
binary trees, where
$n=\dim\di\Ga$
is the topological dimension of the boundary at infinity of
$\Ga$.
\end{abstract}

\section{Introduction}

A metric tree is a geodesic metric space in which every
triangle is isometric to a (maybe degenerate) tripod , i.e.
to the union of three segments with a common point which
is only common point for any pair of the segments.
A metric tree
$T$
is simplicial if it admits a triangulation.
In this case, we can speak about vertices and edges of
$T$.
We always assume that every edge of a
simplicial metric tree has length 1.
The {\em valence} of a vertex is the number of edges
adjacent to it. The {\em binary} metric tree is a
simplicial metric tree in which the valence of every
vertex equals 3.

\begin{thm}\label{thm:main1} Every Gromov hyperbolic group
$\Ga$
admits a quasi-isometric embedding into the product of
$(n+1)$
copies of the binary metric tree where
$n=\dim\di\Ga$
is the topological dimension of the boundary at infinity.
\end{thm}

This result is optimal in the following strong sense:
Any Gromov hyperbolic group
$\Ga$
with
$|\di\Ga|\ge 3$
admits no quasi-isometric embedding into the
$n$-fold
product of any metric trees,
$n=\dim\di\Ga$,
even if the product is
stabilized by any Euclidean factor
$\R^m$, $m\ge 0$.
This is proven in \cite{BS2}.

To give appropriate credit we remark that
the proof of the result relies in an essential
way on the results and methods of the papers
\cite{BL} and \cite{DS}.

Actually in this paper we only consider
$\Ga$
as a metric space and do not use the group structure.
Indeed we will show the following:

\begin{thm} \label{thm:mainthmcapdim} Let
$X'$
be a visual hyperbolic space such that the boundary
$\di X'$
is a doubling metric space. Then
$X'$
admits a quasi-isometric embedding into the product of
$(n+1)$
copies of the binary metric tree, where
$n$
is now the capacity dimension of the boundary at infinity.
\end{thm}

For the definition of a visual hyperbolic space see
sect.~\ref{subsect:hypspace}.
Any Cayley graph of every hyperbolic group is a visual, cocompact
hyperbolic geodesic space and its boundary at infinity is doubling
with respect to any visual metric. Thus we obtain
Theorem~\ref{thm:main1} from Theorem~\ref{thm:mainthmcapdim}
using the following result of \cite{BL}:

\begin{thm}\label{thm:cdimhypspace} The capacity dimension of
the boundary at infinity of every cocompact, hyperbolic
geodesic space
$X$
coincides with the topological dimension,
$\cdim\di X=\dim\di X$.
\qed
\end{thm}

The result should be compared with the Bonk-Schramm embedding
theorem \cite{BoS}, which itself uses the Assouad embedding
result \cite{As}. This embedding
result can be stated
in the following way

\begin{thm}[Bonk-Schramm]\label{thm:bosch} Let
$X'$
be a visual Gromov hyperbolic space such that the boundary
$\di X'$
is a doubling metric space. Then there is a number
$N\in\N$
such that
$X'$
admits a rough similar embedding into the standard
hyperbolic space
$\H^N$.
\end{thm}

The advantage of the Bonk-Schramm embedding is that
the target space is (in contrast to a product of trees)
itself a hyperbolic space and the property of the
embedding map (rough-similarity) is quite strong.
The dimension
$N$
of the target space depends however on the doubling constant of
$\di X'$.

The advantage of our embedding is that the dimension
of the target space is optimal and depends only on the
topological (resp. capacity) dimension of
$\di X'$
and not on the doubling constant of the metric.
This becomes clear in the following examples.

Consider the hyperbolic buildings
$X(p,q)$, $p\ge 5$, $q\ge 2$,
whose apartments are hyperbolic planes with curvature
$-1$,
whose chambers are regular hyperbolic
$p$-gons
with angle
$\pi/2$
and whose link of each vertex is the complete
bipartite graph with
$q+q$
vertices,
studied by Bourdon \cite{Bou}. Indeed there are
infinitely many quasi-isometry classes of these buildings
(distinguished by the conformal dimension of their boundary).
However all of them admit cocompact group actions and
the topological dimension of its boundary is
$1$.
Thus by our result, they all allow quasi-isometric
embeddings into the product of two binary trees.
Note that a product of two binary trees is (the simplest nontrivial)
affine building of rank 2.

If
$X$
and
$Y$
are hyperbolic and
$f:X\to Y$
is a quasiisometric embedding, then the
conformal dimensions satisfy
$\dim_C(\di X) \le \dim_C(\di Y)$
(see \cite[remarques 1.7]{Bou}).
Thus the existence of a quasi-isometric embedding of
$X(p,q)$
into a hyperbolic space
$\H^N$,
implies
$N-1 = \dim_C(\di \H^N)\ge\dim_C (\di X(p,q))$.
The cited paper contains the estimate
$$\dim_C(\di X(p,q))\ge\frac{\log (q-1) + \log p}{2\log p},$$
hence we see that the dimension of the target space
$\H^N$
has to be arbitrarily large as
$q\to\infty$.

We state some consequence of our result
and some possible direction of research:
In dimension theory, there is an important notion of a dimensionally
full-valued space. One of equivalent definitions (due to P.~Alexandrov)
says that a compact space
$X$
is {\em dimensionally full-valued}, if
$\dim(X\times Y)=\dim X+\dim Y$
for every compact space
$Y$.
For example, every 1-dimensional compact space is
dimensionally full-valued, while there are 2-dimensional
compact spaces, the famous Pontryagin surfaces, which are dimensionally
nonfull-valued.
In the paper \cite{BS2} we introduced a dimension invariant,
the {\em hyperbolic dimension}
$\hypdim$
of a metric space. It has the usual properties of a dimension
and is related to the asymptotic dimension
$\asdim$
by
$\hypdim(X)\le\asdim(X)$.
We have
$\hypdim(\R^n)=0$,
but
$\hypdim(\H^n)=n$.
Thus it is a direct consequence of Theorem~\ref{thm:main1}
that for the hyperbolic dimension of the
$n$-fold
product of the binary tree
$T$,
we have
$$\hypdim(T\times\dots\times T)=n=n\hypdim T.$$
This raises the question, is it true that the binary
tree is dimensionally full-valued for the hyperbolic
dimension in the class of all (proper) metric spaces, that is,
$$\hypdim(X\times T)=\hypdim X+1$$
for every (proper) metric space
$X$?
The same question is seemingly open also for the asymptotic dimension.
As far as we know, currently there is no known example
of metric spaces
$X$, $Y$,
which violates the equalities
$\asdim (X\times Y)=\asdim X+\asdim Y$
or
$\hypdim (X\times Y)=\hypdim X+\hypdim Y$
as well as only in rare cases these equalities are known.

\section{Outline of the proof}

In this section we give an outline of the proof of
Theorem \ref{thm:mainthmcapdim}.

\smallskip
\noindent {\bf 1 Step} {\em Hyperbolic approximation:}
\smallskip

The hyperbolic space
$X'$
has the ideal boundary
$Z =\di X'$.
During the proof we will make all constructions and
calculations purely in the space
$Z$.
In particular we will model also the ``interior''
$X'$
as objects in the boundary
$Z=\di X'$.
This realization of the interior in terms of the boundary
can be easily demonstrated in the upper half space model
of the standard hyperbolic space
$\H^{n+1}$.
A point
$(x,r)\in\R^n\times(0,\infty)$
in the upper half space
$\H^{n+1}$
can be viewed as the ball
$B_r(x)\sub\R^n$
in the ideal boundary.

In our situation the boundary
$Z$
is (in contrary to the upper half space model of
$\H^{n+1}$)
compact, but the analogy still works.
Actually to the metric space
$Z$
we will associate a metric graph
$X$
whose vertex set
$V$
consists of a set of balls in
$Z$.
There is in addition a parameter
$0<r<1$
such that the vertex set
$V$
decomposes as
$V=\cup V_k$,
and every element in
$V_k$
is a metric ball
$B(v)$
of radius
$2r^k$
in
$Z$.
We prove the existence of a quasi-isometric embedding
$X'\to X$
and it holds that
$\di X=Z=\di X'$.

Thus roughly speaking we have represented the space
$X'$
as sets of balls in
$Z$
and the radii of the balls are all of the form
$2r^k$
for integers
$k$.

This first step is carried out in section \ref{sect:hypap}

\smallskip
\noindent {\bf  Step 2:} {\em Embedding of
$X$
in a product of infinite valence trees}
\smallskip

Also this step can be visualized in the easiest way in
the upper half space model of
$\H^{n+1}$
with the Euclidean space as its boundary. Consider a family
of subsets of Euclidean space as depictured in the following
Figure \ref{Fi:1}.

\begin{figure}[htbp]
\centering
\includegraphics[width=1.0\columnwidth]{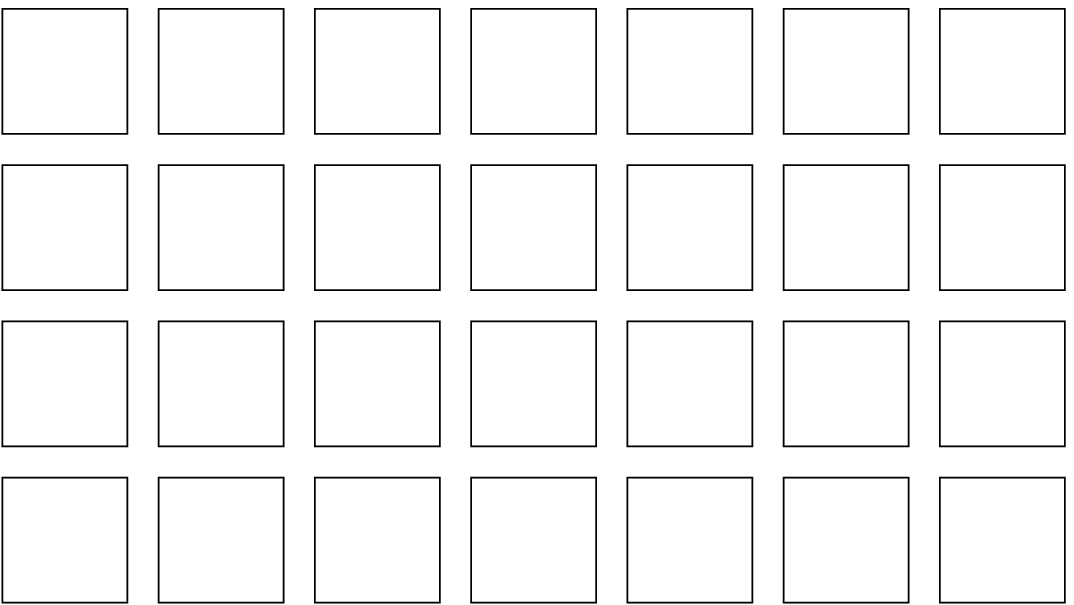}
\caption{Covering sets}\label{Fi:1}
\end{figure}

\begin{figure}[htbp]
\centering
\includegraphics[width=1.0\columnwidth]{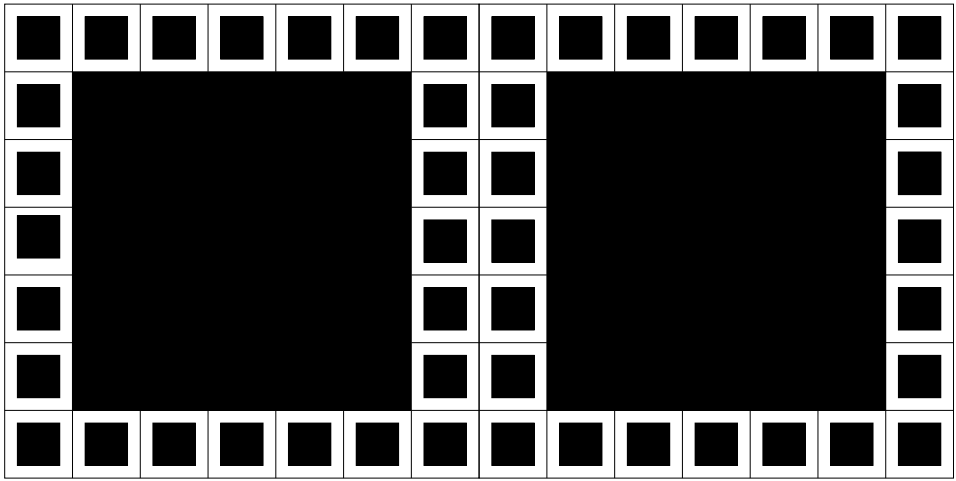}
\caption{Covering sets of neighboring levels}\label{Fi:2}
\end{figure}

It is easy to see that three of these families cover
$\R^2$
or in general
$n+1$
of these families cover
$\R^n$.
In this way one obtains a covering of
$\R^n$
which is colored by
$n+1$
colors such that sets with the same color do not intersect.
Using suitable homotheties one can construct
such families on each scale
$r^k$,
for integers
$k$
and some parameter
$r$
in a way such that for two sets of the same color
(and different scales) either one is contained in
the other or the two sets are disjoint.
More precisely there are
$n+1$
families of subsets
$\cU^c$
of
$X$,
indexed by
$c\in C$, $|C| = n+1$.
For fixed
$c$
the family
$\cU^c$
is a union
$\cU^c=\cup_k\cU^c_k$.
If
$U$, $V\in\cU^c$,
then either
$U\cap V=\es$
or one of the two sets is contained in the other.
The Figure~\ref{Fi:2} shows for one color
$c$
(black in the picture) two neighboring levels
$k$
and
$k+1$.
Note that some smaller cubes of the level
$k+1$
are hidden behind the large black cubes of level
$k$.

Therefore
$\cU^c$
has combinatorially the structure of a tree
(with infinite valence). A vertex
$U\in\cU^c$
can be considered as a large black cube of level
$k$
in Figure~\ref{Fi:2}. Replace this cube by a ``clear window''.
Then draw all black cubes of level
$k+1$, $k+2,\ldots$
which fit into this window. We then see the well picture
of the complement to a Sierpinski carpet formed out of many smaller and
smaller black cubes. All these small cubes that we can see
(i.e. which are not hidden by other black cubes) are
the neighbors with level
$>k$
of
$U$.
We denote this tree by
$T_c$.

If we on the other hand fix
$k$,
then
$\cup_c \cU^c_k$
is an open covering such that the Lebesgue number of this
covering as well as the diameter of each of its elements
is (up to fixed multiplicative constants) of the same size
$r^k$,
i.e. of the same size as the radii of the balls
$B(v)$
for
$v\in V_k$.

It is now possible to define a map
$X\to T_c$.
The map is defined on the set
$V$
of vertices. A vertex
$v$,
which is just the ball
$B(v)$,
is mapped to the smallest
$U\in\cU^c$
such that
$B(v)\sub U$.

It turns out that the product map
$V\to\prod_c T_c$
is a quasi-isometry.

For the standard hyperbolic space this construction
is carried out in \cite{BS1}. For more general spaces
$Z$,
this construction is carried out in \cite{Bu}, where then
$n$
is the capacity dimension of
$Z$.
We discuss the capacity dimension in section \ref{sect:capdim}
and the quasi-isometric embedding into the product of trees in
section \ref{sect:trees} using a simplified version of the
construction from \cite{Bu}. A significant simplification is
achieved due to replacing the hyperbolic
cone construction over
$Z$
used in \cite{Bu} by a hyperbolic approximation of
$Z$
which is much more suitable for this and many other purposes.

\smallskip

\noindent {\bf Step 3:} {\em Alice diary}
\smallskip

This part is completely independent from the rest
of the paper. The main ideas of this construction are
from \cite{DS}. We construct a map from a certain
infinite valence tree into a finite valence tree.
Start with a finite alphabet
$A$
and consider the set
$W$
of words in the alphabet
$A$.
Then
$W$
is an infinite set and hence the corresponding tree
$T_W$
(see section \ref{sect:trees}) of finite sequences in
$W$
is of infinite valence. Alice diary is a certain map
$\psi:T_W\to T_{\Om}$,
where
$\Om$
is finite, hence the tree
$T_{\Om}$
is of bounded valence (and thus quasi-isometric to the binary tree).
In this map the Morse-Thue sequence plays an important role in order
to ``synchronize'' diaries.

This step is contained in the sections \ref{sect:alice} and
\ref{sect:mtsequence}.

\smallskip
\noindent {\bf Step 4:} {\em Labelling of the trees
$T_c$
and the proof of the main theorem}
\smallskip

The tree
$T_c$
is of infinite valence, i.e. a vertex in the tree has
infinite many neighbors. However the neighbors have
different ``levels''. Let us consider a vertex in the tree
$T_c$.
This vertex is just a set
$U\in\cU_k^c$.
A neighbor of level
$k+1$
is a set
$V\in\cU_{k+1}^c$
with
$V \sub U$.
Using the fact that the space
$Z$
is doubling, one can see that there is a constant
$D$
such that
$U$
has at most
$D$
neighbors of level
$k+1$
and more generally only
$D^q$
neighbors of level
$k+q$.
We will more generally show that it is possible to ``label''
the tree
$T_c$
by a finite alphabet. This means that we can isometrically
embed the tree
$T_c$
in a tree
$T_W$,
which is the tree of sentences of a finite alphabet
$A$.
This makes it possible to apply the construction of step 3.
By restricting the diary map
$\psi:T_W\to T_{\Om}$
to the trees
$T_c$
for all colors
$c$,
we obtain the map
$V\to\prod_c T_c\to\prod T_{\Om}$.

It is important to mention that Alice diary map
$\psi:T_W\to T_{\Om}$
is by no means a quasi-isometry. However, the composition
$V\to\prod T_{\Om}$
turns out to be quasi-isometric.

This final step is carried out in section~\ref{sect:colorhypap}.

It remains to note that
$T_\Om$
is quasi-isometric to a subtree of the binary tree,
see Lemma~\ref{lem:binary}
(the tree
$T_\Om$
is a simplicial metric tree whose vertices except the
root have one and the same finite valence
$n\ge 3$.
It can be isometrically embedded into a homogeneous
simplicial metric tree with valence
$n$
of the vertices. On the other hand, it is proven
in \cite{Pa} that every two homogeneous simplicial
metric trees with finite valence
$\ge 3$
are even bilipschitz equivalent to each other).


\section{Preliminaries} \label{sect:prelim}

\subsection{Metric spaces}

Let
$X$
be a metric space. Given points
$x$, $y\in X$
we denote by
$|xy|$
the distance of these points.

A subset
$A\sub Y$
in a metric space
$Y$
is called a {\em net}, if the distances of all points
$y\in Y$
to
$A$
are uniformly bounded.

A map
$f:X\to Y$
between metric spaces is said to be {\em quasi-isometric,}
or a {\em quasi-isometric embedding}
if there are
$a_1$, $a_2 >0$, $b\ge 0$,
such that
$$a_1|xx'|-b\le|f(x)f(x')|\le a_2|xx'|+b$$
for all
$x$, $x'\in X$.
If in addition, the image
$f(X)$
is a net in
$Y$,
then
$f$
is called a {\em quasi-isometry}, and the spaces
$X$
and
$Y$
are called quasi-isometric. In other words, a map is
quasi-isometric if it is bilipschitz on large scales. If
$a_1=a_2$
then
$f$
is called {\em rough similar}.

A metric space
$Z$
is {\em doubling}, if there exists a constant
$C\ge 1$,
such that every ball of radius
$r$
in
$Z$
can be covered by at most $C$
balls of radius
$r/2$.

\subsection{Hyperbolic metric spaces}
\label{subsect:hypspace}

For
$x$, $y$, $z\in X$
we define the {\em Gromov product}
$$(x|y)_z:=\frac{1}{2}(|zx|+|zy|-|xy|).$$
Let
$\de\ge 0$.
A triple
$(a_1,a_2,a_3)\in\R^3$
is called a {\em
$\de$-triple}, if
$$a_{\mu}\ge\min\{a_{\mu+1},a_{\mu+2}\}-\de$$
for
$\mu =1,2,3$,
where the indices are taken modulo 3.

The space
$X$
is called
{\em $\de$-hyperbolic}
if for every
$o$, $x$, $y$, $z\in X$
\begin{equation} \label{eq:de-eq1}
((x|y)_o,(y|z)_o,(x|z)_o)
\ \ \mbox{is a}
\ \de\mbox{-triple}.
\end{equation}

$X$
is called {\em hyperbolic}, if it is
$\de$-hyperbolic
for some
$\de\ge 0$.
The relation~(\ref{eq:de-eq1}) is called the
{\em $\de$-inequality}
with respect to the point
$o\in X$.

If
$X$ satisfies the
$\de$-inequality
for one individual base point
$o\in X$,
then it satisfies the
$2\de$-inequality
for any other base point
$o' \in X$,
see for example \cite{G}. Thus, to check hyperbolicity,
one has to check this inequality only at one point.

Let
$X$
be a hyperbolic space and
$o\in X$
be a base point. A sequence of points
$\{x_i\}\sub X$
{\em converges to infinity,}
if
$$\lim_{i,j\to\infty}(x_i|x_j)_o=\infty.$$
Two sequences
$\{x_i\}$, $\{x_i'\}$
that converge to infinity are {\em equivalent} if
$$\lim_{i\to\infty}(x_i|x_i')_o=\infty.$$
Using the
$\de$-inequality,
one easily sees that this defines an equivalence relation
for sequences in
$X$
converging to infinity. The {\em boundary at infinity}
$\di X$
of
$X$
is defined as the set of equivalence classes of sequences
converging to infinity.

For points
$\xi$, $\xi'\in\di X$
we define their Gromov product by
$$(\xi|\xi')_o=\inf\liminf_{i\to\infty}(x_i|x_i')_o,$$
where the infimum is taken over all sequences
$\{x_i\}\in\xi$, $\{x_i'\}\in\xi'$.
Note that
$(\xi|\xi')_o$
takes values in
$[0,\infty]$
and that
$(\xi|\xi')_o=\infty$
if and only if
$\xi=\xi'$.

If
$\xi$, $\xi'$, $\xi''\in \di X$,
then
$((\xi,\xi')_o,(\xi',\xi'')_o,(\xi,\xi'')_o)$
is a
$\de$-triple.
A metric
$d$
on the boundary at infinity
$\di X$
of
$X$
is said to be {\em visual}, if there are
$o\in X$, $a>1$
and positive constants
$c_1$, $c_2$,
such that
$$c_1a^{-(\xi|\xi')_o}\le d(\xi,\xi')\le c_2a^{-(\xi|\xi')_o}$$
for all
$\xi$, $\xi'\in\di X$.
In this case, we say that
$d$
is the visual metric w.r.t. the base point
$o$
and the parameter
$a$.

The following is well known, see e.g. \cite{BoS}, \cite{V}.

\begin{lem}\label{lem:visualmetric} Let
$X$
be a hyperbolic space. Then for any
$o\in X$,
there is
$a_0>1$
such that for every
$a\in(1,a_0]$
there exists a metric
$d$
on
$\di X$,
which is visual w.r.t.
$o$
and
$a$.
\qed
\end{lem}

A hyperbolic space
$Y$
is said to be {\em visual}, if for some base point
$o\in Y$
there is a positive constant
$D$
such that for every
$y\in Y$
there is
$\xi\in\di Y$
with
$|oy|\le(y|\xi)_o+D$
(one easily sees that this property is independent of the choice of
$o$).
For hyperbolic geodesic spaces this property is a rough version
of the property that every segment
$oy\sub Y$
can be extended to a geodesic ray beyond the end point
$y$.


\section{Hyperbolic approximation of a metric space}
\label{sect:hypap}

Let
$X'$
be a hyperbolic space. Then its boundary
$Z=\di X'$
is a complete metric spaces, since
$Z$
carries a visual metrics with respect to
some base point
$o'\in X'$
and with respect to some parameter
$a'>1$.

In this chapter we start on the other hand
with an arbitrary complete bounded metric space
$(Z,d)$
and construct a geodesic hyperbolic space
$X$
out of
$Z$,
such that
$Z$
can be identified with the boundary at infinity of
$X$.
The space
$X$
is a metric graph and its vertices are balls in
$Z$.
We call
$X$
a hyperbolic approximation of
$Z$.

Finally we prove the following.
Let
$X'$
be a visual hyperbolic space and let
$Z$
be its boundary at infinity.
Let
$X$
be a hyperbolic approximation of
$Z$,
then
$\di X=Z=\di X'$
and there exists a quasi-isometric embedding
$X'\to X$.

\subsection{Construction}\label{subsect:construction}

The construction of a hyperbolic approximation of a metric space
$Z$
is a further development of constructions from \cite{El}
in the case
$Z$
is a compact subspace of an Euclidean space and from \cite{BP}
for arbitrary compact spaces.

Our construction differs from that of \cite{BP} by the definition
of radial edges and radii of balls, which provides some
technical advantages. A hyperbolic approximation can be
defined and turns out to be useful in many situations
for arbitrary metric spaces. For simplicity, we consider
here hyperbolic approximations only of bounded spaces.
Theorem~\ref{thm:approximate} below is similar to
\cite[Proposition 2.1]{BP}.

A subset
$V$
of a metric space
$Z$
is called
$a$-{\em separated}, $a>0$,
if the distance
$d(v,v')\ge a$
for each distinct
$v$, $v'\in V$.
Note that if
$V$
is maximal with this property, then the union
$\cup_{v\in V}B_a(v)$
of balls of radius
$a$
centered at
$v\in V$
covers
$Z$.

Assume that a metric space
$Z$
is bounded,
$\diam Z<\infty$,
and nontrivial, i.e. it contains at least two points.
A {\em hyperbolic approximation} of
$Z$
is a graph
$X$
which is defined as follows. We fix a positive
$r\le 1/6$
which is called the {\em parameter} of
$X$.
Then the largest integer
$k$
with
$\diam Z<r^k$
exists, and we denote it by
$k_0=k_0(\diam Z,r)$.
Note that if
$r<\min\{\diam Z,1/\diam Z\}$
then
$k_0=0$
(the case
$\diam Z<1$)
or
$k_0=-1$
(the case
$\diam Z\ge 1$).

For every
$k\in \Z$, $k\ge k_0$,
let
$V_k\sub Z$
be a maximal
$r^k$-separated
net. One associates with every
$v\in V_k$
the ball
$B(v)\sub Z$
of radius
$r(v):=2r^k$
centered at
$v$.
Note that
$V_{k_0}$
consists of one point. We call this
point the {\em root} of
$X$.

We consider the set
of balls
$B(v)$ for
$v\in V=\cup_{k\ge k_0}V_k$,
as the vertex set of a graph
$X$.

Vertices
$v$, $v'\in V$
are connected by an edge if and only if they either belong to the
same level,
$V_k$,
and the closed balls
$\ov B(v)$, $\ov B(v')$
intersect,
$\ov B(v)\cap\ov B(v')\neq\es$,
or they lie on neighboring levels
$V_k$, $V_{k+1}$
and the ball of the upper level,
$V_{k+1}$,
is contained in the ball of the lower level,
$V_k$.

An edge
$vv'\sub X$
is called {\em horizontal,} if its vertices
belong to the same level,
$v$, $v'\in V_k$
for some
$k\ge k_0+1$.
Other edges are called {\em radial.}
We consider the path metric on
$X$,
for which every edge has length 1. We denote by
$|vv'|$
the distance between point
$v$, $v'\in V$
in
$X$,
and by
$d(v,v')$
the distance between them in
$Z$.
The level function
$\ell :V\to\Z$
is defined by
$\ell(v)=k$
for any
$v\in V_k$.

We often use the following

\begin{rem}\label{rem:net} For every
$z\in Z$
and every
$j\ge k_0$,
there is a vertex
$v\in V_j$
with
$d(z,v)\le r^j$.
This follows from the fact that
$V_j$
is a maximal
$r^j$-separated
set in
$Z$.
\end{rem}

\subsection{Geodesics in a hyperbolic approximation}
\label{subsect:geohypap}
Note that any (finite or infinite) sequence
$\{v_k\}\sub V$,
for which
$v_kv_{k+1}$
is a radial edge for every
$k$,
and the level function
$\ell$
is monotone along
$\{v_k\}$,
is the vertex sequence of a geodesic in
$X$.
Such a geodesic is called {\em radial}.

\begin{lem}\label{lem:radedge} For every
$v\in V$ (except the root vertex)
there is a vertex
$w\in V$
with
$\ell(w)=\ell(v)-1$
connected with any
$v'\in V$, $\ell(v')=\ell(v)$, $|vv'|\le 1$,
by a radial edge. Furthermore,
$d(v,w)\le r^k$
where
$k=\ell(w)$.
\end{lem}

We call the vertex
$w$
a {\em central ancestor} of
$v$.
In general, a central ancestor of
$v$
may not be unique.

\begin{proof} Assume
$v\in V_{k+1}$.
By Remark~\ref{rem:net}, there is a vertex
$w\in V_k$,
for which the distance in
$Z$
between
$v$
and
$w$
is at most
$r^k$, $d(v,w)\le r^k$.
Thus for every vertex
$v'\in V_{k+1}$
adjacent to
$v$
in
$X$,
we have
$$d(v',w)\le d(v',v)+d(v,w)<4r^{k+1}+r^k.$$
For each
$z\in B(v')$
we have
$$d(z,w)\le d(z,v')+d(v',w)<6r^{k+1}+r^k\le 2r^k,$$
since
$r\le 1/6$.
Hence
$B(v')\sub B(w)$,
and
$wv'$
is a radial edge.
\end{proof}

\begin{lem}\label{lem:geoconnect} For every
$v$, $v'\in V$
there exists
$w\in V$
with
$\ell(w)\le\ell(v),\ell(v')$
such that
$v$, $v'$
can be connected to
$w$
by radial geodesics. In particular, the space
$X$
is geodesic.
\end{lem}

\begin{proof} Let
$\ell(v)=k$
and
$\ell(v')=k'$.
Choose
$m<\min\{k,k'\}$
small enough such that
$d(v,v')\le r^{m+1}$.
Applying Lemma~\ref{lem:radedge} we find radial geodesics
$\ga=v_kv_{k-1}\dots v_m$
and
$\ga'=v'_{k'}v'_{k'-1}\ldots v_m'$
in
$X$
connecting
$v=v_k$
and
$v'=v_{k'}'$
respectively with
$m$-th
level. It follows from the definition of
radial edges that
$v\in B(u)$, $v'\in B(u')$
for every vertex
$u\in\ga$, $u'\in\ga'$.
Then
$$d(v',v_m)\le d(v',v)+d(v,v_{m+1})+d(v_{m+1},v_m)
           \le 3r^{m+1}+r^m\le 2r^m$$
since
$r\le 1/6$.
Thus
$v'\in\ov B(v_m)\cap\ov B(v_m')$,
and the vertices
$v_m$, $v_m'$
are connected by a horizontal edge. Applying
Lemma~\ref{lem:radedge} once again we find
$w\in V_{m-1}$
connected with
$v_m$, $v_m'$
by radial edges. Therefore,
$v$, $v'$
are connected to
$w$
by radial geodesics, and
$X$
is connected. This implies that
$X$
is geodesic, because distances between vertices take integer values.
\end{proof}

\begin{lem}\label{lem:horedge} Assume that
$|vv'|\le 1$
for vertices
$v$, $v'$
of one and the same level,
$\ell(v)=\ell(v')$.
Then
$|ww'|\le 1$
for any vertices
$w$, $w'$
adjacent to
$v$, $v'$
respectively and sitting one level below.
\end{lem}

\begin{proof} The balls
$\ov B(w)$, $\ov B(w')$
intersect since they contain the balls
$\ov B(v)$, $\ov B(v')$
respectively, which intersect.
\end{proof}

\begin{lem} Any two vertices
$v$, $v' \in V$
can be joined by a geodesic
$\ga=v_0\ldots v_{n+1}$
such that
$\ell(v_i)<\max\{\ell(v_{i-1}),\ell(v_{i+1})\}$
for all
$1\le i\le n$.
\end{lem}

\begin{proof} Let
$n=|vv'|-1$.
Consider a geodesic
$\ga =v_0\ldots v_{n+1}$
from
$v_0=v$
to
$v_{n+1}=v'$
such that
$\sigma(\ga)=\sum_{i=1}^{n}\ell(v_i)$
is minimal. We claim that
$\ga$
has the desired properties. Let
$1\le i\le n$,
and let
$k=\ell(v_i)$.
Consider the sequence
$(\ell(v_{i-1}),\ell(v_i),\ell(v_{i+1}))$.
There are nine combinatorial possibilities
for this sequence. To prove the result it
remains to show, that the sequences
$(k-1,k,k-1)$,
$(k,k,k)$,
$(k-1,k,k)$
and
$(k,k,k-1)$
cannot occur.

If the sequence is
$(k-1,k,k-1)$,
then
$B(v_i) \sub B(v_{i-1} \cap B(v_{i+1})$ and hence
$|v_{i-1}v_{i+1}|\le 1$
in contradiction to the fact that
$\ga$ is a geodesic.
In the case
$(k,k,k)$ Lemma~\ref{lem:radedge}
implies the existence of
$w\in V_{k-1}$ with
$|v_{i-1}w|\le 1$ and
$|v_{i+1}w|\le 1$.
Replacing the string
$v_{i-1}v_iv_{i+1}$ by
$v_{i-1}wv_{i+1}$
we obtain a new geodesic
$\ga'$
between
$v$, $v'$
with
$\sigma(\ga')<\sigma(\ga)$
in contradiction to the choice of
$\ga$.
The two last cases are symmetric and
we consider only the case
$(k-1,k,k)$.
Choose similar as above
$w\in V_{k-1}$ with
$|v_{i+1}w|\le 1$.
Then by Lemma~\ref{lem:horedge}
$|v_{i-1}w|\le 1$.
Again
$v_{i-1}wv_{i+1}$ defines a geodesic
with smaller
$\si$.
\end{proof}

From this we easily obtain the following

\begin{lem}\label{lem:geohorizone} Any vertices
$v$, $v'\in V$
can be connected in
$X$
by a geodesic which contains at most one horizontal edge.
If there is such an edge, then it lies on the lowest
level of the geodesic.
\qed
\end{lem}

The following corollary is useful in many circumstances.

\begin{cor}\label{cor:ballsintersect} Assume that for some
$v$, $v'\in V$
the balls
$B(v)$, $B(v')$
intersect. Then
$|vv'|\le|\ell(v)-\ell(v')|+1$.
\end{cor}

\begin{proof} We can assume that
$\ell(v)\ge\ell(v')$.
For every vertex
$w\in V$
of a radial geodesic descending from
$v$
we have
$B(v)\sub B(w)$,
in particular, if
$\ell(w)=\ell(v')$
then
$|wv'|\le 1$.
It follows that
$v'$
is the lowest vertex of a geodesic
$v'v\sub X$
as in Lemma~\ref{lem:geohorizone}, hence the claim.
\end{proof}

\subsection{Boundary of the hyperbolic approximation}

The following theorem is a version of \cite[Proposition~2.1]{BP}
adapted to our definition of a hyperbolic approximation,
thus we omit the proof.

\begin{thm}\label{thm:approximate}
Given a complete bounded metric space
$Z$,
its hyperbolic approximation
$X$
is a Gromov hyperbolic geodesic space
with boundary at infinity
$\di X=Z$,
and the metric of
$Z$
is a visual metric on
$\di X$.
The last means that for each
$\xi$, $\xi'\in\di X=Z$
we have
$$c_1a^{-(\xi|\xi')}\le d(\xi,\xi')
  \le c_2a^{-(\xi|\xi')},$$
where positive constants
$c_1$, $c_2$
depend only on
$a= 1/r$
and
$\diam Z$.
\qed
\end{thm}

\subsection{A quasi-isometric embedding}

In this section the following notation is useful.
We write
$A\doteq_C B$
instead of
$|A-B|\le C$.
Sometimes we just write
$A\doteq B$
to indicate that
$|A-B|$
is bounded by some (not specified) constant
$C$.

Let us now assume that
$X'$
is some visual Gromov hyperbolic space with boundary
$Z=\di X'$.
Fix a point
$o'\in X'$.
Let
$d$
be a visual metric on
$Z$
with respect to a parameter
$a'>1$.
Thus there are two positive constants
$c_1$, $c_2$
with
\[c_1 a'^{-(\xi_1'|\xi_2')_{o'}}\le d(\xi_1',\xi_2')\le
 c_2 a'^{-(\xi_1'|\xi_2')_{o'}}\tag{$\ast$}\]
for all
$\xi_1'$, $\xi_2'\in Z$.

Note that if we scale the metric on
$X'$
by some factor
$\la>0$,
all expressions of the form
$(\xi_1'|\xi_2')_{0'}$
are also multiplied by
$\la$.
Thus after a suitable scaling of the metric on
$X'$,
we can assume that
$a'=a$
in estimates
$(\ast)$,
where
$a=1/r\ge 6$
is the constant of the construction in
section~\ref{subsect:construction} above. Note that
the scaling of the metric does not change
the quasi-isometry class of
$X'$.

Let
$X$
be a hyperbolic approximation of the metric space
$(Z,d)$
as constructed above. Then we can identify
$\di X=Z=\di X'$.
We write this identification formally as a map
$f:\di X' \to \di X$.
Using estimates
$(\ast)$
with
$a'=a$
and Theorem~\ref{thm:approximate}, we obtain
$$(f(\xi'_1)|f(\xi'_2))_{o}\doteq_C(\xi'_1|\xi'_2)_{o'}$$
for all
$\xi_1'$, $\xi_2'\in\di X'$
and for some constant
$C$
which depends on
$X$, $X'$
and
$r$,
but not on the arguments
$\xi_i'$.

A map
$F:X'\to X$
between metric spaces is said to be {\em roughly isometric} if
$|F(x)F(x')|\doteq_b |xx'|$
for some constant
$b\ge 0$
and for all
$x$, $x'\in X$.

\begin{thm}\label{thm:roughisometric} Let
$X'$
be a visual and
$X$
a geodesic hyperbolic spaces. Assume that there is a map
$f:\di X'\to\di X$
such that
$(f(\xi'_1)|f(\xi'_2))_{o}\doteq(\xi'_1|\xi'_2)_{o'}$
for all
$\xi_1'$, $\xi_2'\in\di X'$.
Then there exists a roughly isometric map
$F:X\to X'$.
\end{thm}

In the following arguments, all Gromov products
are taken with respect to the base points
$o\in X$
resp.
$o'\in X'$.
To simplify the notation, we omit the base point and denote
$|x|=|ox|$, $|x'|=|o'x'|$
for
$x\in X$, $x'\in X'$
respectively.

The idea of the proof is easily explained in the case
that every point in
$X$
resp.
$X'$
lies on rays emanating from the origin. In this case for
$x'\in X'$
let
$\xi'\in\di X'$ such that
$x'\in o'\xi'$.
Let
$\xi=f(\xi')$
and choose a ray
$o\xi\sub X$.
Then define
$x=F(x')\in o\xi$
to be the point with
$|x|=|x'|$.
Since the equality
$$(\xi'_1|\xi'_2)\doteq(\xi_1|\xi_2)$$
holds up to a uniformly bounded error, one can check that
$F$
is roughly isometric. Under the more general assumptions
of the theorem one has to modify the argument.

We need the following lemma from \cite[Lemma~5.1]{BoS}.

\begin{lem}\label{lem:grprvisual} Let
$X$
be a Gromov hyperbolic space satisfying the
$\de$-in\-equa\-li\-ty
w.r.t. the base point
$o\in X$.
Assume that
$|x_i|\le(x_i|z_i)+D$
for some
$D\ge 0$, $x_i\in X$, $z_i\in X\cup\di X$,
$i=1,2$.
Then
$$(x_1|x_2)\doteq\min\{|\,x_1|,(z_1|\,z_2),|\,x_2|\}$$
up to an error
$\le D+2\de$.
\qed
\end{lem}

\begin{proof}[Proof of Theorem~\ref{thm:roughisometric}]
By the assumption there exists some
$D>0$
such that for every point
$x'\in X'$
there is a point
$\xi'=\xi'(x')\in\di X'$,
with
$(x'|\xi')\ge |x'|-D$.
Choose such
$\xi'$
and let
$\xi=f(\xi')\in\di X$.
Choose a point
$z \in X$
with
$(z|\xi)\ge|x'|$,
in particular,
$|z|\ge|x'|$.
Let
$oz$
be a geodesic from
$o$
to
$z$.
Then we define
$x =F(x')\in oz$
to be the point with
$|x|=|x'|$.
Note that
$F(o')=o$.

Recall that
$$(f(\xi'_1)|f(\xi'_2))\doteq(\xi'_1|\xi'_2)$$
up to an uniform error
$\le c$
for all
$\xi'_1$, $\xi'_2\in\di X'$.

Now given
$x'_i\in X'$
consider
$\xi'_i=\xi'(x'_i)\in\di X'$, $z_i\in X$
with
$(z_i|f(\xi'_i))\ge|x'_i|$
and
$x_i=F(x'_i)\in oz_i$, $i=1,2$.
By Lemma~\ref{lem:grprvisual} we have
$$(x_1'|x_2')\doteq\min\{|x_1'|,(\xi_1'|\xi_2'),|x_2'|\}$$
up to an error
$\le D+2\de'$.
Since
$|x_i'|=|x_i|=(x_i|z_i)$,
we obtain
$$(x_i|f(\xi_i'))\ge\min\{(x_i|z_i),(z_i|f(\xi_i'))\}-\de
  =|x_i|-\de.$$
Then again by Lemma~\ref{lem:grprvisual} we have
$(x_1|x_2)\doteq\min\{|x_1|,(f(\xi_1')|f(\xi_2')),|x_2|\}$
up to an error
$3\de$.
This implies
$(x_1|x_2)\doteq(x_1'|x_2')$
and hence
$|x_1x_2|\doteq|x_1'x_2'|$
up to an error
$\le c+D+2\de'+3\de$.

This shows that
$F$
is roughly isometric.
\end{proof}


\section{Capacity dimension}
\label{sect:capdim}

Let
$Z$
be a metric space. For
$U$, $U'\sub Z$
we denote by
$\dist(U,U')$
the distance between
$U$
and
$U'$,
$\dist(U,U')=\inf\set{d(u,u')}{$u\in U,\ u'\in U'$}$
where
$d(u,u')$
is the distance between
$u$, $u'$.
For
$r>0$
we denote by
$B_r(U)$
the open
$r$-neighborhood
of
$U$, $B_r(U)=\set{z\in Z}{$\dist(z,U)<r$}$,
and by
$\ov B_r(U)$
the closed
$r$-neighborhood
of
$U$, $\ov B_r(U)=\set{z\in Z}{$\dist(z,U)\le r$}$.

Given a family
$\cU$
of subsets in a metric space
$Z$
we define
$\mesh(\cU)=\sup\set{\diam U}{$U\in\cU$}$.
The {\em multiplicity} of
$\cU$, $m(\cU)$,
is the maximal number of members of
$\cU$
with nonempty intersection.
We say that a family
$\cU$
is {\em disjoint} if
$m(\cU)=1$.

A family
$\cU$
is called a {\em covering} of
$Z$
if
$\cup\set{U}{$U\in\cU$}=Z$.
A covering
$\cU$
is said to be {\em colored} if it is the union
of
$m\ge 1$
disjoint families,
$\cU=\cup_{c\in C}\cU^c$, $|C|=m$.
In this case we also say that
$\cU$
is
$m$-colored.
Clearly, the multiplicity of a
$m$-colored
covering is at most
$m$.

Let
$\cU$
be an open covering of a metric space
$Z$.
Given
$z\in Z$,
we let
$L'(\cU,z)=\sup\set{\dist(z,Z\sm U)}{$U\in\cU$}$,
$$L(\cU,z)=\min\{L'(\cU,z),\mesh(\cU)\}$$
be the Lebesgue number of
$\cU$
at
$z$
(the auxiliary
$L'(\cU,z)$
might be larger than
$\mesh(\cU)$
and even infinite as e.g. in the case
$Z=U$
for some member
$U\in\cU$),
$L(\cU)=\inf_{z\in Z}L(\cU,z)$
be the Lebesgue number of
$\cU$.
We have
$L(\cU)\le L(\cU,z)\le\mesh(\cU)$
and for every
$z\in Z$
the open ball
$B_r(z)$
of radius
$r=L(\cU)$
centered at
$z$
is contained in some member of the covering
$\cU$.

There are several equivalent definitions of the capacity
dimension. In this paper, we shall use the
following one. The {\em capacity dimension} of a metric space
$Z$, $\cdim Z$,
is the minimal integer
$m\ge 0$
with the following property: There is a constant
$\de>0$
such that for every sufficiently small
$\tau>0$
there exists a
$(m+1)$-colored
open covering
$\cU$
of
$Z$
with
$\mesh(\cU)\le\tau$
and
$L(\cU)\ge\de\tau$.

\begin{thm}\label{thm:charseq} Let
$Z$
be a bounded metric space with finite capacity dimension,
$n=\cdim Z<\infty$.

Then for every sufficiently small
$r$, $r\in(0,r_0)$,
there exists a sequence
$\cU_j$, $j\ge 0$,
of
$(n+1)$-colored
(by a set
$C$)
open coverings of
$Z$
such that for any hyperbolic approximation
$X$
of
$Z$
with parameter
$r$
the following holds
\begin{itemize}
\item[(1)] $\cU_0^c=\{Z\}$
for all $c\in C$
and
$\mesh\cU_j<r^j$
for every
$j\in\N$;
\item[(2)] for every
$v\in V_{j+1}$, $j\ge 0$,
there is
$U\in\cU_j$
such that
$B(v)\sub U$;

\item[(3)]  for every
$c\in C$
and for different members
$U\in\cU_j^c$, $U'\in\cU_{j\,'}^c$
with
$j\,'\le j$
the following holds: let
$B(U)=\cup\set{B(v)}{$v\in V_{j+1},\ B(v)\cap U\neq\es$}$;
then either
$B(U)\sub U'$
or
$B(U)\cap U'=\es$.
\end{itemize}
\end{thm}

Property (3) is most important and we call it the
{\em separation property}.

\begin{proof} It is proved in \cite[Proposition~2.3]{Bu},
that under the condition of the Theorem, there are constants
$\de$, $\ga\in(0,1)$
such that for every sufficiently small
$r>0$
there exists a sequence
$\cU_j$, $j\in\N$,
of
$(n+1)$-colored
(by a set
$C$)
open coverings of
$Z$
with the following properties
\begin{itemize}
\item[(i)] $\mesh\cU_j<r^j$
and
$L(\cU_j)\ge\de r^j$
for every
$j\in\N$;
\item[(ii)] for every
$c\in C$
and for different members
$U\in\cU_j^c$, $U'\in\cU_{j\,'}^c$
with
$j\,'\le j$
we have either
$B_s(U)\cap U'=\es$,
or
$B_s(U)\sub U'$
for
$s=\ga r^j$.
\end{itemize}
(In \cite{Bu}, the first property of (i) is formulated
as nonstrict inequality; however, the proof actually yields
the strict inequality).
Choosing
$r>0$
sufficiently small, we can assume that
$k_0(\diam Z,r)\le 0$.
Furthermore, we add to the sequence
$\cU_j$, $j\in\N$,
the member
$\cU_0$
which consists of (copies of)
$Z$
for every color
$c\in C$.
Then, property (1) is satisfied.

Because
$L(\cU_j)\ge\de r^j$,
every ball
$B_\rho(z)\sub Z$
of radius
$\rho\le\de r^j$
is contained in some member
$U\in\cU_j$.
Assuming that
$r<\de/2$,
we obtain property (2).

Finally, assume additionally that
$r<\ga/4$.
Now, if a color
$c\in C$
and different members
$U\in\cU_j^c$, $U'\in\cU_{j\,'}^c$
with
$j\,'\le j$
are given, we have
$B(U)\sub B_s(U)$
with
$s=\ga r^j$
by the choice of
$r$
and the definition of
$B(U)$.
Hence, property (3).
\end{proof}


\section{Trees}
\label{sect:trees}

\subsection{Levelled trees}\label{subsect:leveltree}

Recall a poset (partially ordered set)
$V$
is called {\em directed}, if for any
$u$, $v\in V$
there is
$w\in V$
with
$u\le w$, $v\le w$.

A {\it levelled tree}
$T$
is a directed poset
$V$,
called the {\em vertex set} of
$T$,
together with a {\em level} function
$\ell:V\to\Z$, which is strictly monotone
in the following sense:
If $v, v' \in V$ are different elements
and $v\le v'$, then
$\ell(v)>\ell(v')$.

In this case,
$v'$
is called an {\em ancestor} of
$v$,
and
$v$
is a {\em descendant} of
$v'$.

We require that the following condition is satisfied:

$(+)$
if distinct elements
$v$, $v'\in V$
have a common descendant then one of them is an ancestor
of the other.

A collection
$E$
of two point subsets of
$V$
called the {\em edge set} of
$T$
is defined by the condition: A pair of vertices
$(v,v')$
forms an {\em edge},
$(v,v')\in E$,
if and only if one of its member, say
$v'$,
is an ancestor of the other and the level
$\ell(v')$
is maximal with this property.

If there is a vertex, which has no ancestor,
then such a vertex is unique by directedness, and it is called the
{\em root} of
$T$.
Note that the root is an ancestor of every other vertex.

It follows from
$(+)$
that for every vertex
$v\in V$
(except the root) there is exactly one edge
$(v,v')$,
in which
$v$
is the descendant. Hence, by the uniqueness part,
$T$
has no circuit. By the existence part (together with
$(+)$),
every vertex is connected with any its ancestor
by a sequence of edges in
$T$.

Now, because
$V$
is directed, every two vertices in
$T$
are connected by a sequence of edges, i.e.
$T$
is connected. Therefore,
$T$
is a simplicial tree.

The members of an edge are called
{\em neighbors}. So, the number of neighbors of a vertex
$v$
equals the valence of
$v$.

In an obvious way, one can define a distance
function on the vertex set of a levelled tree.
The distance
$\cL(v,v')$
between
$v\in V$
and its ancestor
$v'\in V$
is the number of {\em generations} between them, i.e.
$\cL(v,v')=k$,
if there is a sequence
$v'=v_0,\dots,v_k=v\in V$
such that
$(v_i,v_{i+1})\in E$
and
$v_i$
is the ancestor of
$v_{i+1}$
for
$i=0,\dots,k-1$.
In particular,
$\cL(v,v')\le\ell(v)-\ell(v')$
and the inequality might be strong. In general case,
the distance
$\cL(v,v')$
between vertices is the sum of their distances
to the youngest common ancestor
$v_0$, $\cL(v,v')=\cL(v,v_0)+\cL(v_0,v')$.
If
$T$
is rooted and
$\es$
is the root, we define
$\cL(v):=\cL(v,\es)$.
Therefore, each levelled tree is a simplicial metric tree
with an interior metric with length 1 edges.

With every alphabet
$\Om$,
one associates the rooted levelled tree
$T_\Om$
in the following way: its vertex set
$V$
is the set of all words in
$\Om$,
i.e. the set of all finite sequences with elements in
$\Om$,
with the obvious partial order:
$\psi\le\psi'$
for words
$\psi$, $\psi'$
if and only if
$\psi'$
is an initial subword of
$\psi$.
This order is obviously directed and condition
$(+)$
is satisfied. The value of the level function
$\ell(\psi)$
is the number of letter in the word
$\psi$.
Then two vertices
$\psi$, $\psi'\in V$
are connected by an edge in
$T_\Om$
if and only if one of them (ancestor) is the initial
string of the other (descendant) obtained by erasing
the last member of the descendant.

The empty sequence defines the root of
$T_\Om$.
The root has
$|\Om|$
neighbors and every other vertex has
$|\Om|+1$
neighbors, one ancestor and
$|\Om|$
descendants. Here,
$|\Om|$
denotes the cardinality of
$\Om$.
Note that
$\ell(\psi)$
coincides with the length of the word
$\psi$, $\ell(\psi)=\cL(\psi)$,
in this case.

\begin{lem}\label{lem:binary} For every finite alphabet
$\Om$,
the tree
$T_\Om$
is roughly homothetic a subset of the binary tree.
\end{lem}

\begin{proof} We can assume that
$\Om=\{1,\dots,n\}$
for some
$n\ge 3$.
Every number
$k\in\Om$
can be uniquely written as the string of zeros and ones
of length
$\la=[\log_2n]+1$
via the binary representation (with an appropriate number
of zeros in front if necessary). This defines a map
$f:T_\Om\to T_{\{0,1\}}$,
where
$T_{\{0,1\}}$
is the {\em rooted} binary tree, which is obviously
isometric to a subtree of the binary tree.

The map
$f$
is radially homothetic with coefficient
$\la$,
$\cL(f(\psi),f(\psi'))=\la\cL(\psi,\psi')$
for each pair (ancestor, descendant). It easily
follows from this that in general case, we have
$$\la\left(\cL(\psi,\psi')-2\right)+2\le
  \cL(f(\psi),f(\psi'))\le\la\cL(\psi,\psi').$$
\end{proof}

For an arbitrarily rooted levelled tree,
we have only
$\ell(v)\le\cL(v)$
for its vertices
$v$.
Moreover, vertices at a fixed distance to the root
may have different and even arbitrarily large
levels. Important examples of
such trees are colored trees associated with
sequence of colored coverings as in Theorem~\ref{thm:charseq}.

\subsection{Colored trees}\label{sect:coltree}

Let
$Z$
be a bounded metric space with finite capacity dimension,
$n=\cdim Z<\infty$;
$X$
be a hyperbolic approximation of
$Z$
with sufficiently small parameter
$r<\min\{\diam Z,1/\diam Z\}$
satisfying the condition of Theorem~\ref{thm:charseq};
$\{\cU_j\}$, $j\ge 0$,
be a sequence of
$(n+1)$-colored
(by a set
$C$)
open coverings of
$Z$
with parameter
$r$
as in Theorem~\ref{thm:charseq}.
Then recall
$k_0=k_0(\diam Z,r)=0$
or
$-1$.

For every
$c\in C$,
we define a rooted levelled tree
$T_c$
as follows. Its vertex set
$\cU^c$
is the disjoint union
$\cU^c=\cup_{j\ge 0}\cU_j^c$
with the root
$v_c=Z$,
which is the unique member of
$\cU_0^c$.
The partial order of
$\cU^c$
is defined by the inclusion relation,
$U\le U'$
if and only if
$U\sub U'$,
and this order is directed due to the existence of the root.

We say that a vertex
$U\in\cU_j^c$
has level
$j$;
this defines the level function.
It follows easily from the separation property (3)
that the level function is strictly monotone.
Then a vertex
$U\in\cU_j^c$
is a descendant of
$U'\in\cU_{j\,'}^c$
if
$j\,'<j$
and
$U\sub U'$.

It follows from the separation property (3) that
$(+)$
is satisfied. Hence,
$T_c$
is a rooted levelled tree: a pair of vertices
$U\in\cU_j^c$, $U'\in\cU_{j\,'}^c$
forms an edge of
$T_c$
if and only if it is a pair (descendant, ancestor),
and the level
$j\,'$
of the ancestor
$U'$
is maximal with this property.

Typically, vertices of the tree
$T_c$
have infinite valence, see e.g. Introduction and \cite{BS1}.

We use notation
$\cL(U,U')$
for the distance in
$T_c$
between its vertices
$U$, $U'$,
and
$\cL(U)$
for the distance
$\cL(U,v_c)$.
Note that
$\cL(v_c)=0$
and
$\cL(U)\le j-k_0$
if
$U\in\cU_j^c$.
Furthermore, if
$U'\in\cU_{j\,'}^c$
is an ancestor of
$U\in\cU_j^c$,
the level difference
$j-j'$
might be arbitrarily large compared to
the distance
$\cL(U,U')$
even if
$(U,U')$
is an edge of
$T_c$.
Keeping this in mind is highly useful for
understanding of what follows.

\subsection{A map into the product of colored trees}
\label{subsect:maprodtree}

We now define a map
$f_c:V\to T_c$,
$V$
is the vertex set of the hyperbolic approximation
$X$
of
$Z$,
as follows. The root of
$X$
is mapped into the root of
$T_c$, $f_c(v)=v_c$
for the unique member
$v\in V_{k_0}$.
Given
$v\in V_j$, $j > k_0$,
we let
$f_c(v)=U\in\cU_{j\,'}^c$
be the covering element containing the ball
$B(v)$, $B(v)\sub U$,
and
$j\,'\le j-1$
is maximal with this property. By Theorem~\ref{thm:charseq}(3),
$f_c(v)$
is well defined.

\begin{lem}\label{lem:lipschitz} For every
$c\in C$,
the map
$f_c:V\to T_c$
is Lipschitz,
$$\cL(f_c(v),f_c(v'))\le 2|vv'|,\quad
  \text{for every}\quad v, v'\in V.$$
\end{lem}

\begin{proof} Since the hyperbolic approximation
$X$
is geodesic, it suffices to estimate the distance
$\cL(f_c(v),f_c(v'))$
for neighbors,
$|vv'|=1$.

Assume that the edge
$(v,v')\sub X$
is horizontal, i.e.
$v$, $v'\in V_j$
for some
$j\ge k_0$
and the balls
$B(v)$, $B(v')$
intersect. Thus the covering elements
$U=f_c(v)$, $U'=f_c(v')$
also intersect. By the separation property (3),
either
$U=U'$
and so
$f_c(v)=f_c(v')$
or these elements have different levels and one of
them is contained in the other, say
$U\in\cU_i^c$, $U'\in\cU_{i'}^c$
with
$i>i'$
and
$U\sub U'$.
It follows from the definition of
$f_c$
that
$i<j$
and from the separation property (3) that any
$U''\in\cU_{i''}^a$, $i''<i$,
intersecting
$U$,
also contains
$B(v')$.
Thus
$(U,U')\sub T_c$
is an edge by the definition of
$f_c$,
and
$\cL(f_c(v),f_c(v'))=1$
in this case.

Assume now that the edge
$(v,v')\sub X$
is radial, say
$v\in V_{j+1}$, $v'\in V_j$.
Then
$B(v)\sub B(v')$
and as in the previous case,
$U\cap U'\neq\es$.
Thus we can assume that these elements have
different levels,
$U\in\cU_i^c$, $U'\in\cU_{i'}^c$
with
$i\neq i'$,
and one of them is contained in the other. Moreover,
$i>i'$
because
$B(v)\sub U'$,
thus
$U\sub U'$.
We have
$i\le j$
by definition of
$f_c$,
and as above, any ancestor
$U''\in\cU_{i''}^c$, $i''<i-1$,
of
$U$
separated from
$U$
by at least one generation,
$\cL(U,U'')\ge 2$,
also contains
$B(v')$.
Therefore, it follows from the definition of
$f_c$
that at most one generation can separate
$U$
from its ancestor
$U'$,
and
$\cL(f_c(v),f_c(v'))\le 2$.
\end{proof}

\begin{thm}\label{thm:firstqiso}
The map
$$f=\prod_c f_c:V\to\prod_c T_c$$
is quasi-isometric.
\end{thm}

The map
$f:V\to\prod_{c\in C}T_c$
defined by its coordinate maps
$f_c:X\to T_c$
is Lipschitz by Lemma~\ref{lem:lipschitz}.
To prove that
$f$
is roughly bilipschitz, we begin with the following Lemma,
which is the main ingredient of the proof.

For
$i\ge 0$,
we denote by
$T_{c,i}=\cU_i^c$
the vertex set of
$T_c$
of level
$i$.

\begin{lem}\label{lem:radial} Given
$v\in V_{j+1}$, $j\ge 0$,
for every integer
$i$, $0\le i\le j$,
there is a color
$c\in C$
such that
$\cL(f_c(v),T_{c,i})\ge M$
with
$M+1\ge(j-i+1)/|C|$.
Furthermore, if for
$k\le i$
a vertex
$w\in T_{c,k}$
is the lowest level vertex of the segment
$f_c(v)w\sub T_c$,
then
$\cL(f_c(v),w)\ge M$.
\end{lem}

\begin{proof} Consider a radial geodesic
$v_{i+1}\dots v_{j+1}\sub X$
with vertices
$v_m\in V_m$,
where
$v_{j+1}=v$.
This means in particular that
$B(v_{m+1})\sub B(v_m)$
for every
$m=i+1,\dots,j$.
By Theorem~\ref{thm:charseq}(2) for every vertex
$v_{m+1}$,
there is a covering element
$U_m\in\cU_m$
with
$B(v_{m+1})\sub U_m$, $m=i,\dots,j$.

There is a color
$c\in C$
such that the set
$\{U_i,\dots,U_j\}$
contains
$M+1\ge(j-i+1)/|C|$
members having the color
$c$,
i.e. every of those
$U_m\in\cU_m^c$.
Since
$B(v)\sub U_m$
for every
$m\le j$,
we have
$U=f_c(v)\sub U_m$
for every
$U_m$
having the color
$c$
by the definition of
$f_c$
and the separation property (3).

Using again the separation property, we obtain
that any path in
$T_c$
between
$f_c(v)$
and the set
$T_{c,i}$
must contain at least
$M+1$
vertices and hence
$\cL(f_c(v),T_{c,i})\ge M$.

Finally, let
$W\in\cU_k^c$
be the set corresponding to the vertex
$w$
of
$T_c$.
By the assumption on
$w$,
the set
$W$
contains
$U$
and every set from the list
$\{U_i,\dots,U_j\}$
having the color
$c$.
Hence,
$\cL(f_c(v),w)\ge M$.
\end{proof}

We say that distinct points
$v\in V_j$, $v'\in V_{j\,'}$, $j\ge j\,'\ge 0$
are {\em horizontally close} to each other, if
$d(v,v')<r^{j\,'}$.
This term is motivated by the fact that there is
a geodesic segment
$vv'\sub X$
which is almost radial. More precisely, we have

\begin{lem}\label{lem:radclose} Assume that the distinct points
$v$, $v'\in V$
are horizontally close to each other. Then their levels are
different and the upper level ball is contained in the
lower level ball, say
$\ell(v)>\ell(v')$, $B(v)\sub B(v')$.
In particular,
$|vv'|\le|\ell(v)-\ell(v')|+1.$
\end{lem}

\begin{proof} We can assume that
$v\in V_j$, $v'\in V_{j\,'}$, $j\ge j\,'\ge 0$.
Then
$j>j\,'$
because
$v$, $v'$
are distinct and because
$V_j\sub Z$
is
$r^j$-separated
for every
$j\ge 0$.
Furthermore,
$B(v)\sub B(v')$
because
$2r^j+r^{j\,'}<2r^{j\,'}$.
By Corollary~\ref{cor:ballsintersect}, we have
$|vv'|\le(j-j\,')+1$.
\end{proof}

\begin{pro}\label{pro:estbelowhorclose} Given
$v$, $v'\in V$
horizontally close to each other, we have:
$f_c(v)f_c(v')\sub T_c$
is a radial segment for every color
$c\in C$,
i.e. its lowest level vertex is one of its ends,
and there is a color
$c\in C$
such that
$$|vv'|\le|C|\cL(f_c(v),f_c(v'))+\si,$$
where
$\si=|C|+1$.
\end{pro}

\begin{proof} We can assume that
$v$, $v'$
are distinct. Then by Lemma~\ref{lem:radclose},
their levels are different, say
$\ell(v)>\ell(v')$,
and
$B(v)\sub B(v')$.
Thus
$f_c(v)=f_c(v')$
or
$f_c(v)$
is a descendant of
$f_c(v')$
for every color
$c\in C$.
In any case,
$f_c(v')$
is the lowest level vertex of the segment
$f_c(v)f_c(v')\sub T_c$.

On the other hand,
$f_c(v')\in T_{c,m_c}$
with
$m_c\le\ell(v')$
for all
$c\in C$
by the definition of
$f_c$.
By Lemma~\ref{lem:radial}, we have
$$\cL(f_c(v),f_c(v'))=\cL(f_c(v),T_{c,m_c})\ge M$$
with
$M+1\ge(\ell(v)-\ell(v'))/|C|$
for some color
$c\in C$.
Therefore, using again Lemma~\ref{lem:radclose},
we obtain
$$|vv'|\le|\ell(v)-\ell(v')|+1\le
 |C|(M+1)+1\le|C|\cL(f_c(v),f_c(v'))+\si.$$
\end{proof}

\subsubsection{Digression: Critical level of two vertices}
\label{subsubsect:critlevel}

We say that vertices
$v$, $v'\in V$
are {\em horizontally distinct}, if
$\ell(v)$, $\ell(v')\ge 0$
and they are not horizontally close to each other,
$d(v,v')\ge r^{\min\{\ell(v),\ell(v')\}}$.
In this case, there is an integer
$l$
with
$$r^l\le d(v,v')<r^{l-1}.$$
We call
$l$
the {\em critical level} of
$v$, $v'$.
Note that
$r^l\le\diam Z<r^{k_0}$
and thus
$l>k_0$,
in particular,
$l\ge 0$.
Furthermore,
$l\le\min\{\ell(v),\ell(v')\}$.

\begin{lem}\label{lem:critleveldist} Let
$l$
be the critical level of horizontally distinct vertices
$v$, $v'\in V$.
Then
$|vv'|\le\ell(v)+\ell(v')-2l+3$.
\end{lem}

\begin{proof} Consider central ancestor radial geodesic
$\ga_v$
in
$X$
between the root of
$X$
and
$v$.
For the vertex
$u\in\ga_v$
of the level
$l-1$, $u\in V_{l-1}$,
we have
$|vu|=\ell(v)-l+1$
and
$d(v,u)\le r^{l-1}$.
Thus
$d(u,v')
 \le d(u,v)+d(v,v')<2r^{l-1}$.
Therefore,
$v'\in B(u)$,
the balls
$B(u)$, $B(v')$
intersect and hence,
$|v'u|\le\ell(v')-l+2$
by Corollary~\ref{cor:ballsintersect}. We see that
$|vv'|\le |vu|+|uv'|\le\ell(v)+\ell(v')-2l(v,v')+3$.
\end{proof}

\begin{lem}\label{lem:critlevel} Let
$l$
be the critical level of horizontally distinct
$v$, $v'\in V$.
Assume that
$v\in U$, $v'\in U'$
for some elements
$U$, $U'\in T_c$
and some color
$c\in C$.
Then, for the lowest level vertex
$W\in T_{c,k}$
of the segment
$UU'\sub T_c$,
we have
\begin{itemize}
\item[(1)] its level
$k<l$;
\item[(2)] there are at most three vertices on each
of the segments
$WU$, $WU'\sub UU'$
(including
$W$)
having the level
$<l$.
\end{itemize}
\end{lem}

\begin{proof} (1) The set
$W\in\cU_k^c$
contains both
$U$
and
$U'$,
and we have by Theorem~\ref{thm:charseq}(1)
$$d(v,v')\le\diam W<r^k.$$
It follows
$r^l<r^k$
and thus
$k<l$.

(2) Assume there are distinct vertices
$W'$, $W''\in WU$
with levels
$k''<k'<l-1$,
$W'\in T_{c,k'}$, $W''\in T_{c,k''}$.
By the separation property (3),
$B(W')\sub W''$,
where
$B(W')$
is the union of all balls
$B(u)$
intersecting
$W'$
with
$u\in V_{k'+1}$.
Since
$k'+1<l$,
the radius of every such ball is
$2\rho\ge 2r^{l-1}$.
There is
$u\in V_{k'+1}$
with
$d(u,v)\le\rho$.
Then,
$$d(u,v')\le d(u,v)+d(v,v')<\rho+r^{l-1}\le 2\rho,$$
hence
$v'\in B(u)\sub W''$.
It follows
$U'\sub W''$
and thus
$W''=W$.
\end{proof}

The first part of this Lemma will be used in the
proof of Proposition~\ref{pro:estbelow}; the second one
in the proof of Proposition~\ref{pro:crlevelsentences}.

\subsubsection{Proof of Theorem~\ref{thm:firstqiso}}
\label{subsubsect:proofthm}

The following Proposition is a key step in the proof of
Theorem~\ref{thm:firstqiso} and Theorem~\ref{thm:compose}.
Note that its condition is asymmetric with respect to
the points
$v$, $v'$.

\begin{pro}\label{pro:estbelow} Given horizontally distinct
$v$, $v'\in V$, $\ell(v)\ge\ell(v')$,
there is a color
$c\in C$
such that the following two facts hold

\begin{itemize}
\item[(1)] $\max\{\ell(f_c(v)),\ell(f_c(v'))\}-l+1
           \le |C|(\cL(f_c(v),w)+1)$;
\item[(2)] $|vv'|\le 2|C|\cL(f_c(v),w)+\si$,
\end{itemize}
where
$l=l(v,v')$
is the critical level of
$v$, $v'$;
$w\in f_c(v)f_c(v')\sub T_c$
is the lowest level vertex;
$\si=2|C|+1$.
\end{pro}

\begin{proof} By Lemma~\ref{lem:critlevel}(1) for every color
$c\in C$,
any path in
$T_c$
between
$f_c(v)$
and
$f_c(v')$
passes through a vertex of a level
$k<l$.

By Lemma~\ref{lem:radial}, there is a color
$c\in C$
such that
$\cL(f_c(v),T_{c,l-1})\ge M$
with
$M+1\ge(\ell(v)-l+1)/|C|$.
Let
$w\in T_{c,k}$
be the lowest level vertex of the segment
$f_c(v)f_c(v')$.
By Lemma~\ref{lem:radial}, we have
$\cL(f_c(v),w)\ge M$
since
$k\le l-1$.
Because
$\ell(f_c(v))\le\ell(v)-1$
and
$\ell(f_c(v'))\le\ell(v')-1\le\ell(v)-1$
by the definition of the map
$f_c$,
we have
$$\max\{\ell(f_c(v)),\ell(f_c(v'))\}-l+1
   \le\ell(v)-l\le|C|(M+1),$$
hence (1). Because
$\ell(v)\ge\ell(v')$,
we have by Lemma~\ref{lem:critleveldist},
$$|vv'|\le 2(\ell(v)-l+1)+1\le 2|C|\cL(f_c(v),w)+\si,$$
hence (2).
\end{proof}

We use also the notation
$|ww'|$
for the distance between
$w$, $w'\in\prod_cT_c$
(on a product of trees we take the $l_1$ product metric).
The following Proposition completes the proof
of Theorem~\ref{thm:firstqiso}.

\begin{pro}\label{pro:bilip} There are constants
$\La>0$, $\si\ge 0$
depending only on
$|C|$
such that
$$|vv'|\le\La|f(v)f(v')|+\si$$
for all
$v$, $v'\in V$.
\end{pro}

\begin{proof} If vertices
$v$, $v'$
are horizontally close to each other, then the required estimate
follows from Proposition~\ref{pro:estbelowhorclose}.
If they are horizontally distinct, then the required estimate
follows from Proposition~\ref{pro:estbelow}(2) because
$\cL(f_c(v),w)\le\cL(f_c(v),f_c(v'))$.
\end{proof}

\section{Alice Diary} \label{sect:alice}

In this section we construct a map from some
infinite valence tree into some finite valence tree.

\subsection{Sentences in an alphabet}
\label{subsect:sentences}
Start with a finite set
$A$,
which we consider as some alphabet. By
$W$,
we denote the set of all finite words in
$A$.
In particular, the empty word
$\es$
is in
$W$.

Thus
$W$
is the vertex set of
$T_A$,
where we recall from section~\ref{subsect:leveltree}
the definition of the tree
$T_{\Om}$
for some alphabet
$\Om$.

A {\em sentence}
$\al$
in the alphabet
$A$
is a word in the alphabet
$W$.
It is convenient to introduce a special marker
$s$
called the {\em stop sign}, which terminates words
in a sentence. So, we write
$$\al=w_1sw_2s\dots sw_ks$$
for a sentence of
$k$
words. Let
$S$
be the set of all sentences. Thus
$S$
is the vertex set of the tree
$T_W$.

Certainly, the number of letters from the alphabet
$A$
in a sentence
$\al$
may be much larger than the number of words.
Furthermore, the set
$W$
considered as the alphabet for sentences is infinite.
Now, we want to have a method which allows to
encode sentences by a finite alphabet in a way that
any encoded sentence would have exactly the same
number of letters as the initial one has words.
In other words we want to find a finite alphabet
$\Om$
and a map
$\psi:T_W \to T_{\Om}$
preserving the combinatorial distance to the roots.

This is of course impossible if we would try to
retain all of the information contained in every
sentence. However, if we pursue a moderate aim only
to retain the information contained in every
sentence having a definite positive percentage
of stop signs with respect to the number of
letters from
$A$,
the task becomes solvable irrespective of the length
of the words in a sentence.

\subsection{Alice diary}
\label{subsect:diary}
The solution is the Alice diary. Consider a sentence
$\al$
as description of a journey of Alice who
wants to write a diary about her trip. Every letter of
$\al$
represents a day. There are two types of days during
this journey. The days when the weather is so fine that
Alice has no time to write her diary. These are all unmarked
days. Then there are the days of rest marked by
$s$.
In the morning of every
$s$-day,
Alice writes a page describing
$\ka$
days of her journey: she starts with yesterday and then the day before
yesterday etc. Of course she skips
the days which were already described earlier in the diary.
If there is no day left to describe, she marks on the page the symbol $\star$
and stop to write the diary at this rest day.

More formally, fix
$\ka\in\N$
called the {\em diary constant} and consider the alphabet
$\Om'$,
which consists of two types of letters: either
$\om\in\Om'$
is a word in the alphabet
$A'=A\cup\{s\}$
having precisely
$\ka$
letters, or
$\om$
consists of
$<\ka$
letters from
$A'$
concluded by the symbol
$\star$.
We call letters of
$\Om'$
{\em pages}.

Now, given a sentence
$\al$
in the alphabet
$A$
consisting of
$k$
words,
$\al=w_1sw_2s\dots sw_ks$,
we define inductively the diary
$\psi(\al)$
consisting of
$k$
pages in
$\Om'$
as follows. Its first page
$\om_1$
consist of
$\ka$
letters of the word
$w_1$
written backward starting from the last
letter of
$w_1$;
in the case
$\ell(w_1)<\ka$
(i.e. the length of the word is smaller than $\ka$),
the page
$\om_1$
is the word
$w_1$
written backward and augmented by the symbol
$\star$.

Next, we delete
$\om_1$
from
$\al$,
proceed in the same way starting from the last letter of
$w_2$
and obtain the second page
$\om_2$
of
$\psi(\al)$.
After
$k$
steps, the diary
$\psi(\al)$
is completed.

As an example consider
$\ka = 3$
and the sentence
$$\al=\un{a\,a\,b\,c}\,s\,
     \un{a}\,s\,
     \un{b\,c\,b}\,s\,
     \un{c}\,s\,
     \un{b}\,s$$
(words are underlined). Then
$\om_1=(c\,b\,a)$,
$\om_2=(a\,s\,a)$,
$\om_3=(b\,c\,b)$,
$\om_4=(c\,s\,s)$,
$\om_5=(b\,s\,\star)$
and
$$\psi(\al)=\om_1\om_2\om_3\om_4\om_5
  =(c\,b\,a)(a\,s\,a)(b\,c\,b)(c\,s\,s)(b\,s\,\star).$$

Thus we have constructed a map
$T_W\to T_{\Om'}$ respecting the distance to
the roots in the corresponding trees.


\subsection{Reconstruction procedure}
\label{subsect:reconstruct}
Certainly, different sentences in
$A$
may have one and the same diary. Thus, one can reconstruct
a sentence out of its diary only up to this ambiguity.
Given a diary
$\tau$
considered as a word in
$\Om'$,
we let
$A_{\tau}$
be the set of all sentences in
$A$
with the diary
$\tau$,
i.e.
$A_{\tau}=\psi^{-1}(\tau)$.
Note that
$A_{\tau}$
may be empty.

We first explain the direct reconstruction
using the example above, where the last diary entry
$\om_5$
contains the letter
$\star$.
The first entry
$\om_1$
of the diary shows that the subword
$w_1$
has the form
$$\sqcup\,a\,b\,c$$
where the
``$\sqcup$''
stand for a sequence of letters of unknown
length (maybe length 0). The second entry implies that
$w_1sw_2$
has the form
$$\sqcup\, a\,a\,b\,c\,s\,\un{a};$$
the information
$\om_3$
then gives that
$w_1sw_2sw_3$
is of the form
$$\sqcup\,a\,a\,b\,c\,s\,\un{a}\,s\,\sqcup\,b\,c\,b\,;$$
the entry
$\om_4$
says that
$w_1sw_2sw_3sw_4$
has the form
$$\sqcup\,a\,a\,b\,c\,s\,\un{a}\,s\,\un{b\,c\,b}\,s\,\un{c}$$
and the last entry implies that
$w_1sw_2sw_3sw_4sw_5$
has the form
$$\un{a\,a\,b\,c}\,s\,\un{a}\,s\,
  \un{b\,c\,b}\,s\,\un{c}\,s\,\un{b}.$$
Thus we have reconstructed the whole sentence
(we can finally add the stop symbol at the end).
The reconstructed words appear underlined only at
the moment when further steps could not change them.
Note that in this direct recovery procedure, it may
happen that up to the last step no word would be
underlined even if finally the whole sentence is reconstructed.

\begin{rem}\label{rem:eveningdiary}
There is a subtlety in the definition of the
diary. It is important that Alice writes her diary in
the morning and not in the evening. If Alice would write
the diary in the evening and first describe
the rest day, then the reconstruction is not possible:
e.g. the following journeys
$\un{a\,b}\,s\,\un{c\,d}\,s\,\un{e}\,s$
and
$\un{c\,a\,b}\,s\,\un{d}\,s\,\un{e}\,s$
have the same evening-diary
$(s\,b\,a)$, $(s\,d\,c)$, $(s\,e\,\star)$.
\end{rem}

The above reconstruction procedure motivates the
following terminology. We call an expression
$\sqcup\,w$,
where
$w$
is a word, a {\em slotted word}.
A {\em slotted sentence} is an expression of the form
$\be=u_1su_2s\ldots u_ks$,
where now every
$u_i$
is either a word or a slotted word. To every
slotted sentence
$\be$,
we associate the set
$A_{\be}\sub S$
of sentences obtained from
$\be$
by putting any words in the place of the slots.

For example, if
$$\be=\sqcup\,w_1sw_2s\sqcup\,w_3s,$$
where the
$w_i$
are some words, then
$A_{\be}$
consists of all sentences of the form
$$w'_1w_1sw_2s w'_3 w_3s,$$
where
$w'_1$
and
$w'_3$
are arbitrary words.

Clearly
$A_{\be} =\{\be\}$
in the case that
$\be$
is a honest sentence (without slots), and
$A_{\be}$
can be canonically identified with
$W\times\ldots\times W$ ($p$
factors), if
$\be$
has
$p\ge 1$
slots. The
$i$-th
$W$
factor corresponds to the word which is put into the
$i$-th
slot. To be more formal, let
$A_{\be}\to W\times\ldots\times W$
be this identification. We consider also the map
$W\times\ldots\times W\to S$
defined by
$(w_1,\ldots,w_p)\mapsto w_1sw_2s\ldots w_ps$.
Let
$r_{\be}:A_{\be}\to S$
be the composition of these two maps. In the case
that
$\be$
is a honest sentence, let
$r_{\be}(\be)=s$.

In other words, every
$\al\in A_\be$
is obtained from the slotted sentence
$\be$
by filling in its slots by the appropriate words of
the sentence
$r_\be(\al)$.

To formulate the reconstruction lemma, we recall
the algorithm of Alice diary. In each step
Alice has given a string of letters (ending with some
$s$).
In this step, she writes the
$\ka$
last (or
$<\ka$
letters and
$\star$)
letters before the final stop sign in her diary.
We call the string which is left (including the final stop sign)
the {\em rest sentence} of the step.
In our example the rest sentence in the first step is
$a\,s$.
If
$\al$
is a sentence, we denote by
$r(\al)$
the rest sentence after writing the whole diary.
Note that
$r(\al)$
always contains the final
$s$
of
$\al$,
and if the last page of the diary
$\psi(\al)$
is concluded by
$\star$
then whole sentence
$\al$
except the final
$s$
is on the diary and thus
$r(\al)=s$.

\begin{lem}\label{lem:slot} For every sentence
$\al\in S$,
there exists a well defined slotted sentence
$\wh\al$
such that
$$\psi^{-1}(\psi(\al))=A_{\wh\al}\quad\text{and}\quad
r(\al')=r_{\wh\al}(\al')$$
for
every
$\al'\in A_{\wh\al}$.
\end{lem}

\begin{proof} Given a sentence
$\al=w_1\,s\dots w_k\,s$, $k\ge 1$,
we let
$\tau=\psi(\al)$
be the diary of
$\al$, $\tau=\om_1\dots\om_k$.
We inductively construct the required slotted sentence
$\wh\al$.
On every step
$i=1,\dots,k$,
the slotted sentence
$\wh\al_i$
and the rest sentence
$r_i$
are defined using data of the previous step.
We put
$\al_i=w_1\,s\dots w_i\,s$.
It turns out that
$r(\al_i)=r_{\wh\al_i}(\al_i)$
for every
$i=1,\dots,k$,
and we use notation
$r_i$
for this sentence. Then
$\wh\al:=\wh\al_k$
and
$r(\al)=r_{\wh\al}(\al)$.

\smallskip\noindent
1st step. We distinguish two cases:

\noindent {\bf Case 1}:
$w_1$
has
$\ge\ka$
letters,
$w_1=u_1'u_1$,
where
$u_1$
consists of
$\ka$
letters and
$u_1'$
might be empty. We put
$$\wh\al_1=\sqcup\,u_1\,s.$$
The sentence
$\al_1=w_1\,s$
is obtained from the slotted sentence
$\wh\al_1$
by putting its initial string
$u_1'$
into the slot, thus
$\al_1\in A_{\wh\al_1}$.
We have
$$r_{\wh\al_1}(\al_1)=u_1's=r(\al_1),$$
and we denote this sentence by
$r_1$.
In this case,
$\om_1$
is
$u_1$
written in reverse order.

\noindent {\bf Case 2}:
$w_1$
has less than
$\ka$
letters. We put
$$\wh\al_1=w_1\,s.$$
Then
$\wh\al_1$
is a honest sentence, thus
$A_{\wh\al_1}=\{\wh\al_1\}$
and
$\al_1\in A_{\wh\al_1}$.
Then again
$$r_1:=r_{\wh\al_1}(\al_1)=s=r(\al_1).$$
In this case,
$\om_1$
is
$w_1$
written in reverse order followed by
$\star$.

\smallskip\noindent
$(i+1)$-th
step. Assume that we have already defined the slotted
sentence
$\wh\al_i$
such that the sentence
$\al_i=w_1\,s\dots w_i\,s$
is in the class
$A_{\wh\al_i}$, $\al_i\in A_{\wh\al_i}$,
we also have defined the rest sentence
$r_i=r_{\wh\al_i}(\al_i)$,
which coincides with
$r(\al_i)$, $r_i=r(\al_i)$,
and have written the diary pages
$\om_1\dots\om_i$, $i\ge 1$,
where each page depends only on the slotted sentence
of the same step.

We study the sentence
$r_iw_{i+1}s$,
which plays the basic role in the induction step.
We distinguish two cases:

\noindent {\bf Case 1}:
$w_{i+1}$
has at least
$\ka$
letters,
$w_{i+1}=u_{i+1}'u_{i+1}$,
where
$u_{i+1}$
is a word of
$\ka$
letters and the word
$u_{i+1}'$
might be empty.
We put
$$\wh\al_{i+1}=\wh\al_i\sqcup u_{i+1}s.$$
Then using the induction assumption, we see that
$\al_{i+1}\in A_{\wh\al_{i+1}}$,
and we define
$$r_{i+1}:=r_{\wh\al_{i+1}}(\al_{i+1})=r_i\,u_{i+1}'s.$$
Using that
$r_i=r(\al_i)$,
we obviously have
$r(\al_{i+1})=r_i\,u_{i+1}'s$
and thus
$r(\al_{i+1})=r_{i+1}$.
In this case, the page
$\om_{i+1}$
is
$u_{i+1}$
written in reverse order.

\noindent {\bf Case 2}:
$w_{i+1}$
has less than
$\ka$
letters. We distinguish two subcases.

\noindent Subcase 2a:
$\wh\al_i$
is a honest sentence. In particular,
$\al_i=\wh\al_i$
and
$r_i=r_{\wh\al_i}(\al_i)=s=r(\al_i)$.
Thus
$r_iw_{i+1}=sw_{i+1}=:\wh\om_{i+1}$
consists of
$\le\ka$
letters. We put
$$\wh\al_{i+1}=\wh\al_i\,w_{i+1}s.$$
Then
$\wh\al_{i+1}=\al_{i+1}$
is a honest sentence, and we see again that
$\al_{i+1}\in A_{\wh\al_{i+1}}$.
Now, we define
$$r_{i+1}:=r_{\wh\al_{i+1}}(\al_{i+1})=s.$$
The rest sentence
$r(\al_{i+1})$
is obtained from the basic sentence
$r_iw_{i+1}s$
by removing the string
$\wh\om_{i+1}$,
thus
$r(\al_{i+1})=s=r_{i+1}$.
In this case, the page
$\om_{i+1}$
is
$\wh\om_{i+1}$
written in reverse order and, in the case
$\wh\om_{i+1}$
consists of less than
$\ka$
letters, followed by
$\star$.

The most complicated is

\noindent Subcase 2b: the slotted sentence
$\wh\al_i$
contains at least one slot.

Let
$\wh\om_{i+1}$
be the string of the last
$\ka$
letters of
$r_iw_{i+1}$
in the case
$r_iw_{i+1}$
contains at least
$\ka$
letters or otherwise
$\wh\om_{i+1}=r_iw_{i+1}$.
Using
$r_i=r_{\wh\al_i}(\al_i)$,
we see that
$$\wh\om_{i+1}=v_p\,s\dots v_1\,s\,w_{i+1}$$
for some (may be empty) words
$v_1,\dots,v_p\in W$,
where
$p\ge 1$
does not exceed the number of slots in
$\wh\al_i$.

Now, we define the slotted sentence
$\wh\al_{i+1}$
by adding the word
$w_{i+1}\,s$
at the end of
$\wh\al_i$,
filling in the
$j$-th
slots (from the end) of
$\wh\al_i$
by the corresponding
$v_j$, $1\le j\le p-1$,
i.e. replace the last slot of
$\wh\al_i$
by
$v_1$
etc, and finally filling in the
$p$-th
slot (from the end) of
$\wh\al_i$
by
$\sqcup\,v_p$
in the case
$\wh\om_{i+1}$
consists of
$\ka$
letters or otherwise filling in the
$p$-th
slot by
$v_p$.

By the induction assumption,
$r_i=r_{\wh\al_i}(\al_i)$,
that is,
$\al_i$
is obtained from
$\wh\al_i$
by filling in its slots by the appropriate words of
$r_i$.
Thus it immediately follows from the definition of
$\wh\al_{i+1}$
that
$\al_{i+1}\in A_{\wh\al_{i+1}}$.
Now, we define
$$r_{i+1}:=r_{\wh\al_{i+1}}(\al_{i+1}),$$
which is equal to deleting the substring
$\wh\om_{i+1}$
from
$r_iw_{i+1}s$.
By the same operation, the rest sentence
$r(\al_{i+1})$
is obtained from
$r_iw_{i+1}s$
due to the induction assumption. Hence,
$r(\al_{i+1})=r_{i+1}$.

In this case, the diary page
$\om_{i+1}$
is
$\wh\om_{i+1}$
written in reverse order and, in the case
$\wh\om_{i+1}$
consists of less than
$\ka$
letters, followed by
$\star$.
Thus
$\om_{i+1}$
is defined only by a part of
$\al_{i+1}$
present in the slotted sentence
$\wh\al_{i+1}$.

This completes the induction step and hence the
construction of the slotted sentence
$\wh\al$.
By construction,
$\al\in A_{\wh\al}$, $r(\al)=r_{\wh\al}(\al)$
and moreover
$\psi(\al')=\psi(\al)$
for every
$\al'\in A_{\wh\al}$.

Conversely, assume that
$\psi(\al')=\psi(\al)$
for some sentence
$\al'\in S$,
and consider the slotted sentence
$\wh\al'$
constructed as described for
$\al'$.
The inspection of writing the diary pages for
$\tau=\psi(\al)=\psi(\al')$
on every inductive step shows that
$\wh\al'=\wh\al$,
thus
$\al'\in A_{\wh\al}$
and
$r(\al')=r_{\wh\al}(\al')$.
\end{proof}

Every page
$\om$
of the diary
$\psi(\al)$
for some sentence
$\al$
has been written in the corresponding rest day,
which we denote by
$s_{\om}$.

The very last case of the proof above tells us
that if some page
$\om_i$, $i=1,\dots,k$,
of the diary
$\psi(\al)$
contains
$\star$
as the last symbol, then
$\wh\al_i$
is a honest sentence and hence
$A_{\wh\al_i}=\{\wh\al_i\}$.
Thus we obtain

\begin{lem}\label{lem:completerecover} Assume that some page
$\om$
of the diary
$\psi(\al)$
contains
$\star$,
and let
$\wh\al'$
be the slotted sentence from Lemma~\ref{lem:slot}
constructed by the initial subsentence
$\al'\sub\al$
up to the stop sign
$s_{\om}$
corresponding to
$\om$.
Then
$\wh\al'=\al'$
is a honest sentence, thus
$A_{\wh\al'}=\{\al'\}$
and the sentence
$\al'$
is unambiguously reconstructed out of its diary.
\qed
\end{lem}

From a concrete diary
$\tau = \psi(\al)$
consisting out of
$k$
pages we can construct the slotted sentence
$\wh\al$
with
$k$
words. We call
$\wh\al$
the reconstruction of the sentence
$\al$.
If
$\wh\al$
is a honest sentence, then
$\wh\al=\al$
and we say that we have reconstructed the whole word.
In general the reconstruction
$\wh\al$
can be obtained from
$\al$
by replacing certain substrings of
$\al$
by a slot
$\sqcup$
and inserting slots in some places.

If the frequency of the stop sign
$s$
is high in some region of this sentence, then strings
decoded in the diary come together and it is possible to reconstruct
large substrings of
$\al$.

We denote by
$s_m$
the
$m$-th
stop sign of
$\al$
and by
$\om_m$
the corresponding page of
$\psi(\al)$.
The set of letters between
$s_m$
and
$s_{m+p}$
is denoted by
$(s_m,s_{m+p})$
(neither
$s_m$
nor
$s_{m+p}$
is included). As usual, the sign
$\#$
means cardinality, and
$\ka$
the diary constant.

\begin{lem}\label{lem:stringrecover} Assume that
$\ka p\ge\#(s_m,s_{m+p})+1$
for some sentence
$\al\in S$
and
$m$, $p\in\N$.
Then either all or at least
$k$
letters of
$\al$
in a row left to
$s_{m+1}$
are written in its reconstruction
$\wh\al$
as a string without slots inside, where
$k=\ka+q+\#(s_m,s_{m+1})$, $q=\ka p-\#(s_m,s_{m+p})$.
\end{lem}

\begin{proof} If any of the consecutive pages
$\om_{m+1},\dots,\om_{m+p}$
of the diary
$\psi(\al)$
contains
$\star$
then the statement follows from Lemma~\ref{lem:completerecover}.
Thus, we assume that no page above contains
$\star$.
Then, they together contain
$\ka p$
letters of
$\al$
scattered left to
$s_{m+p}$.
Note that the whole substring of
$\al$
from the leftmost of these scattered letters up to
$s_{m+1}$
is on the diary and moreover, it is written in
$\wh\al$
without slots inside. If
$\ka p\ge\#(s_m,s_{m+p})+1$
then
$s_m$
is on the diary among other
$\ka p$
letters. If the page
$\om_m$
contains
$\star$,
then the proof is completed again by
Lemma~\ref{lem:completerecover}. Otherwise,
$\ka+\ka p-\#(s_m,s_{m+p})+\#(s_m,s_{m+1})=k$
letters of
$\al$
in a row left to
$s_{m+1}$
are written in
$\wh\al$
without slots inside.
\end{proof}

\section{The Morse-Thue sequence and synchronization}
\label{sect:mtsequence}

\begin{exa}\label{exa:exodus} The Exodus:
Consider the journey of two brothers, which travel
40 years through the desert until they reach the
Promised Land.
However one of the brothers stays one week longer in the desert.
Let
$w=b\,a\,a\,a\,a\,a\,a$
represents one week. Now, consider the sentences
$$\al=w^k\,s\,\dots\,s$$
where there are
$n$
entries of
$s$
at the end (the first
$n$
days in the Promised Land);
$$\al'=w^{k+1}\,s\,\dots\,s$$
where there are
again
$n$
entries
$s$
at the end. Note that if
$n\ka\ll k$,
then the diaries of two brothers coincide.
\end{exa}

This example shows that even if a diary encodes an arbitrarily
large string of letters, there is no way in general to recover
the place of the string in a sentence, and two different sentences
may have the same diary encoding mutually shifted strings.
The problem is that every letter in a sentence
$\al$
has its own level which in general is lost in the diary. To match
the letters of
$\al$
by levels and retain this information in the diary would require
again an infinite alphabet. We rectify this problem introducing
so called Morse-Thue decoration.

\subsection{Morse-Thue sequence and decoration of sentences}
\label{subsect:mtdecor}

The following sequence was introduced
and studied by A.~Thue \cite{Th} and later independently by
M.~Morse \cite{Mo}.

\begin{Def}\label{Def:mtsequence} Consider the substitution rule
$0\to 01$
and
$1\to 10$.
Then start from
$0$
to perform this substitutions
$$0\to 01\to 0110\to 01101001\to\dots$$
to obtain a nested family of sequences
$t_k$
of length
$2^k$
in the alphabet
$\{0,1\}$.
The resulting limit sequence is called the {\em Morse-Thue}
sequence.
\end{Def}

The Morse-Thue sequence
$\set{t(n)}{$n\ge 0$}$
has the following remarkable
property (see e.g. \cite{He}), which is only used in
what follows.

\begin{thm}\label{thm:cubefree} The Morse-Thue sequence is
cube-free, i.e. it contains no string of type
$w\,w\,w$
where
$w$
is any word in
$0$
and
$1$.
\qed
\end{thm}

Let
$\al$
be a sentence in the alphabet
$A$.
We define the level of its letters inductively
putting
$\lv(a)=1$
for the first letter if it is not the stop sign
and
$\lv(a)=0$
otherwise. Now assume that the level of a letter
$a$
is defined, and
$a'\in\al$
is the next letter. If
$a'$
is not the stop sign then we put
$\lv(a')=\lv(a)+1$,
otherwise,
$\lv(a')=\lv(a)$.
In other words, the level of each stop sign is the same
as the level of the last letter of the word which it
terminates, while the level of any other letter is
defined in the natural way.

Now, we define the Morse-Thue decoration of a sentence
$\al$
as follows. This is the sentence obtained from
$\al$
replacing its
letter
$a$
of the level
$m$
by the letter
$(a,t(m))$
for all
$m\ge 0$.
The diary
$\psi(\al)$
of a decorated sentence
$\al$
is defined exactly as above via the decorated
stop signs, and Lemma~\ref{lem:stringrecover}
holds for decorated sentences. Decorated pages are
now elements of the finite alphabet
$\Om=\Om'\times\{0,1\}$,
and a decorated diary is a string of decorated
pages.

\subsection{Synchronization}\label{subsect:synchron}

The length of a sentence
$\al$, $|\al|$,
is the level of its last letter. In this sense, the stop
signs are ignored while computing the length. We also
ignore stop signs while saying about length of strings
in a sentence unless the opposite is explicitly stated.

Given a decorated sentence
$\al$
and a letter
$a\in\al$,
we denote
by
$t_a$
the {\em tail} of
$\al$
with initial letter
$a$,
that is all letters of
$\al$
following and including
$a$.

\begin{lem}\label{lem:synchron} Assume that decorated
sentences
$\al$, $\al'$
have identical tails of length
$\ge l$
and their lengths differ at
most by one half of
$l$, $||\al|-|\al'||\le l/2$.
Then
$|\al|=|\al'|$.
\end{lem}

\begin{proof} Comparing the decorations of the identical tails,
we immediately observe a subsequence
$w\,w\,w$
of the Morse-Thue sequence unless the lengths of sentences
coincide and therefore the tails are not shifted with respect
to each other.
\end{proof}

Combining recovery Lemma~\ref{lem:stringrecover}
and Lemma~\ref{lem:synchron}, we obtain the following
proposition, which plays an important role in the proof of our
main theorem.

\begin{pro}\label{pro:equaldiaries} Assume for
some decorated sentences
$\al$, $\al'$
we know that
\begin{itemize}
\item[(1)] the sentences
$\al$, $\al'$
both contain no empty word, i.e. there are no
two stop signs in a row in
$\al$, $\al'$;

\item[(2)] there are letters
$a\in(s_m,s_{m+1})\sub\al$
and
$a'\in(s_{m'},s_{m'+1})\sub\al'$
with
$\lv(a)=\lv(a')$
and
$|m-m'|\le 2$;
\item[(3)] there are at least
$p$
stop signs behind
$a$, $a'$
in
$\al$, $\al'$
respectively and
$\max\{|t_a|,|t_{a'}|\}\le n(p-2)$
for some
$p\ge 3$
and
$n\in\N$;

\item[(4)]
$\psi(\al)=\psi(\al')$
for
$\ka$-diaries
with
$\ka\ge 5n+1$.
\end{itemize}
Then
$a=a'$.
\end{pro}

\begin{proof} We can assume that
$m\le m'$.
First, we check the condition of recovery
Lemma~\ref{lem:stringrecover} for the interval
$(s_{m+2},s_{m+p})\sub\al$
(this choice is motivated by the estimate
$m+3\ge m'+1$).
We have
$$q:=\ka(p-2)-\#(s_{m+2},s_{m+p})
   \ge (5n+1)(p-2)-(|t_a|+p-3)\ge 4n(p-2)>1.$$
Thus by Lemma~\ref{lem:stringrecover}, either all or at least
$k$
letters in a row left to
$s_{m+3}$
in
$\al$
are written in its reconstruction
$\wh\al$
as a string without slots inside, and therefore
this string can be recovered from
$\psi(\al)$,
where
$k=\ka+q+\#(s_{m+2},s_{m+3})$.

Because the diaries coincide,
$\psi(\al)=\psi(\al')$,
the reconstruction of
$\al'$
is
$\wh\al$
by Lemma~\ref{lem:slot} and thus
again the string of at least
$k$
letters in a row left to
$s_{m+3}$
in
$\al'$
coincides with the corresponding string of
$\al$.
Therefore, the corresponding tails of
$\al_{m+3}$, $\al_{m+3}'$
coincide, where
$\al_{m+3}$, $\al_{m+3}'$
are the heads of
$\al$, $\al'$
respectively consisting of
$m+3$
words.

Now we check the condition of the synchronization
Lemma~\ref{lem:synchron}. Since there is no empty
word in
$\al$, $\al'$,
the identical recovered tails of
$\al_{m+3}$, $\al_{m+3}'$
both have length
$l\ge k/2\ge 2n(p-2)$.
On the other hand,
$|t_a|\le n(p-2)\le l/2$.
Thus the recovered string of
$\al$
contains
$a$.
Because
$m\le m'\le m+2$,
the
$(m+3)$-th
stop sign in
$\al$
is identical with
$(m'+1)$-th
or with
$(m'+2)$-th
or with
$(m'+3)$-th
stop sign in
$\al'$.
In any case, the recovered string of
$\al'$
contains
$a'$
because also
$|t_{a'}|\le l/2$.
It follows from (2),
$\lv(a)=\lv(a')$,
that
$$||\al_{m+3}|-|\al_{m+3}'||\le\max\{|t_a|,|t_{a'}|\}\le l/2,$$
and we can apply Lemma~\ref{lem:synchron}. Thus
$|\al_{m+3}|=|\al_{m+3}'|$
and hence,
$a=a'$.
\end{proof}


\section{Labelling of a hyperbolic approximation}
\label{sect:colorhypap}

Now, we assume that our metric space
$Z$
is doubling. Then, we can find a finite set
$F$
with cardinality depending only on the doubling constant of
$Z$
such that for every
$j\ge k_0$
there is a coloring
$\mu_j:V_j\to F$
with
$\mu_j(v)\neq\mu_j(v')$
for each distinct
$v$, $v'\in V_j$
with
$d(v,v')<2r^{j-2}$.
Let
$A'$
be the set consisting of all nonempty subsets of
$F$.
Clearly
$A'$ is finite.
The set
$A=A'\times\{0,1\}$
will serve as an alphabet in a way that edges of every tree
$T_c$, $c\in C$,
will be labelled by words in
$A$
and vertices of
$T_c$
by sentences in
$A$.

Namely, to the edge between two vertices
$U\in\cU_j^c$ and $U'\in\cU_{j\,'}^c$
with
$j\,'<j$,
we associate a word in
$A$
of length
$j-j\,'$.
The word consists of letters
$a_k\in A$
for each level
$k$
with
$j\,'+1\le k\le j$.
We define
$a_k=(a_k',t(k))$,
where
$$a_k'=\set{\mu_{k+1}(v)\in F}{$v\in V_{k+1},\ B(v)\cap U\ne\es$}
      \sub F$$
and
$t(k)$
is the
$k$-th
member of the Morse-Thue sequence. Every
vertex
$U$
of the tree
$T_c$
is connected to the root by a unique sequence of
edges and in this way it becomes a sentence
$\al=\al(U)$
in the alphabet
$A$
(or a Morse-Thue decorated sentence in the alphabet
$A'$).
We also insert stop signs
$s$
in
$\al$
terminating the words. The stop signs are only needed
to define the diary of
$\al$.
So, they play an auxiliary role, and their levels
coincide with levels of the last letters in
corresponding words of
$\al=w_1s\dots w_ms$.
Recall that the number of words
$m$
in
$\al$
does not exceed
$j=\sum_i|w_i|$,
where
$j$
is the level
of
$U$,
and
$j$
might be arbitrarily large compared with
$m$.

In this way, the role of every member
$U\in T_c$
is threefold. First, it is a covering element from
$\cU^c$,
in particular, a subset of
$Z$.
Second, it is a vertex of the tree
$T_c$
and thus one can speak about its properties as a member
of the tree. Finally, it can also be considered as
a sentence in the alphabet
$A$.
Depending of a role, we sometime use different notations
for one and the same member of
$T_c$.

Now, we shall use notations of sect.~\ref{sect:alice}
related to sentences.

\begin{pro}\label{pro:crlevelsentences} Let
$\al$, $\al'$
be sentences in the alphabet
$A$
corresponding to vertices
$U$, $U'\in T_c$
respectively. Assume that there are horizontally distinct
$v$, $v'\in V$,
for which
$v\in U$, $v'\in U'$
(as points and subsets in
$Z$), $\ell(U)$, $\ell(U')\ge l+1$
where
$l=l(v,v')$
is the critical level of
$v$, $v'$.
Then for the letters
$a\in(s_m,s_{m+1})\sub\al$,
$a'\in(s_{m'},s_{m'+1})\sub\al'$
of the critical level,
$\lv(a)=l=\lv(a')$,
we have
\begin{itemize}
\item[(1)] $|m-m'|\le 2$;
\item[(2)] $a\neq a'$.
\end{itemize}
\end{pro}

\begin{proof} The condition
$\ell(U)$, $\ell(U')\ge l+1$
implies that the assumption of the lemma is not void, i.e.
the letters
$a$, $a'$
of
$\al$, $\al'$
respectively, exist. Furthermore,
$\diam U$, $\diam U'<r^{l+1}$
by Theorem~\ref{thm:charseq}(1).

Let
$W\in T_c$
be the lowest level vertex of the segment
$UU'\sub T_c$, $k$
the level of
$W$,
i.e.
$W\in T_{c,k}$.
It follows from Lemma~\ref{lem:critlevel}, that
$k<l$
and that there are at most three vertices (including
$W$)
on each of the segment
$WU$, $WU'\sub UU'$
having the level
$<l$.
Hence,
$|m-m'|\le 2$.

To prove (2), we show that there
exists
$u\in V_{l+1}$
such that
$B(u)\cap U\ne\es$, $B(u)\cap U'=\es$
and
$d(u,u')<2r^{l-1}$
for every
$u'\in V_{l+1}$
with
$B(u')\cap U'\ne\es$.
Then, the color from
$F$
corresponding to
$u$, $\mu_{l+1}(u)$,
which is on the list in
$a$,
will not appear in the color list of
$a'$
by the property of the coloring
$\mu_{l+1}$,
and thus
$a\neq a'$.

There exists
$u\in V_{l+1}$
with
$d(u,v)\le r^{l+1}$.
Then
$v\in B(u)$
and thus
$B(u)\cap U\ne\es$.
On the other hand, for every
$z'\in U'$,
we have
$$d(u,z')\ge d(v,v')-d(u,v)-d(v',z')\ge r^l-2r^{l+1}>2r^{l+1}$$
because
$r\le 1/6$.
Hence,
$B(u)\cap U'=\es$.

Assume now that
$B(u')\cap U'\ne\es$
for some
$u'\in V_{l+1}$.
Then
$d(v',u')\le\diam U'+2r^{l+1}<3r^{l+1}$,
and we have
$$d(u,u')\le d(u,v)+d(v,v')+d(v',u')
   <4r^{l+1}+r^{l-1}<2r^{l-1}.$$
Hence, the claim.
\end{proof}

\subsection{The diary map}\label{diarymap}

We fix a natural
$\ka$,
whose value will be specified below, and let
$\Om=\Om'\times\{0,1\}$
be the set of the Morse-Thue decorated pages
associated with sentences in the alphabet
$A$,
see sect.~\ref{subsect:diary} and \ref{subsect:mtdecor}.
Note that one and the same alphabet
$A$
has been used for labelling of every tree
$T_c$, $c\in C$.
Thus, the {\em finite} alphabet
$\Om$
and therefore the tree
$T_\Om$
is independent of colors
$c\in C$.

Now, for every
$c\in C$,
the diary
writing procedure defines a map
$$\psi_c:T_c\to T_\Om,$$
called the {\em diary} map. It follows from
the definition that
$\ell(\psi_c(\al))=\cL(\al)$
for every sentence
$\al\in T_c$,
i.e. the map
$\psi_c$
is radially isometric. Consequently,
$\psi_c$
is 1-Lipschitz and the product map
$$\psi:\prod_{c\in C}T_c\to T_{\Om}^{|C|}$$
is also Lipschitz (this map is 1-Lipschitz if both products,
$\prod_{c\in C}T_c$
and
$T_{\Om}^{|C|}$,
are equipped with the same type of product metrics, either
$\ell_1$,
$\ell_2$
or
$\ell_\infty$,
what is natural to consider). Certainly,
$\psi$
as well as every
$\psi_c$, $c\in C$,
is by no means quasi-isometric. However, we have

\begin{thm}\label{thm:compose} The composition map
$\psi\circ f:V\to T_{\Om}^{|C|}$
is quasi-isometric.
\end{thm}

\begin{proof} We use notations
$\eta=\psi\circ f$
and
$\eta_c=\psi_c\circ f_c$
for every color
$c\in C$.
The map
$\eta$
is Lipschitz because both
$f$
and
$\psi$
are Lipschitz. To prove that
$\eta$
is roughly bilipschitz, we consider two cases.
The easy case is if
$v$, $v'\in V$
are horizontally close to each other. Then
by Proposition~\ref{pro:estbelowhorclose}, we have:
$f_c(v)f_c(v')\sub T_c$
is a radial segment for every color
$c\in C$,
and there is a color
$c\in C$
such that
$$|vv'|\le|C|\cL(f_c(v),f_c(v'))+\si,$$
where
$\si=|C|+1$.
Because
$\psi_c$
is radially isometric,
$\cL(f_c(v),f_c(v'))=\cL(\eta_c(v),\eta_c(v'))$,
and therefore
$$|vv'|\le|C||\eta(v)\eta(v')|+\si.$$

The case
$v$, $v'$
are horizontally distinct is much more interesting.
We can assume that
$|vv'|\ge 15|C|^2+2|C|+\si$,
where now
$\si=2|C|+1$
according Proposition~\ref{pro:estbelow}(2), and
$\ell(v)\ge\ell(v')$.

Then, by Proposition~\ref{pro:estbelow}, there is a color
$c\in C$
such that
\begin{itemize}
\item[(1)] $\max\{\ell(f_c(v)),\ell(f_c(v'))\}-l+1
          \le|C|(\cL(f_c(v),w)+1)$;
\item[(2)] $\cL(f_c(v),w)\ge(|vv'|-\si)/|C|\ge 15|C|+2$,
\end{itemize}
where
$w\in f_c(v)f_c(v')\sub T_c$
is the lowest level vertex,
$l=l(v,v')$
is the critical level of
$v$, $v'$.
We fix this color
$c$
and consider it in what follows. Furthermore, we use notation
$\be=f_c(v)$, $\be'=f_c(v')$
for elements
$f_c(v)$, $f_c(v')\in T_c$
considered as sentences in the alphabet
$A$.

Using the terminology of
sect.~\ref{sect:alice},
$q=\cL(f_c(v),w)$
is the number of stop signs behind
$w$
in the sentence
$\be$.
By Lemma~\ref{lem:critlevel}(2), only at most
two of them are sitting below the critical level
$l$.
By (2), $q-2\ge 15|C|\ge 15$,
thus there is the letter
$a\in\be$
of the critical level,
$\lv(a)=l$.
Then in notations of sect.~\ref{sect:alice},
$a\in(s_m,s_{m+1})\sub\be$
for some
$m\ge 1$.
Hence, there are at least
$q-2$
stop signs behind
$a$
in
$\be$.

Furthermore, (1) above means that the tail
$t_a(\be)\sub\be$
of
$\be$
with initial letter
$a$
has length
$|t_a(\be)|=\ell(f_c(v))-l+1\le|C|(\cL(f_c(v),w)+1)$.
Because no empty word occurs in
$\be$,
we have
$|t_a(\be)|\ge q-2\ge\cL(f_c(v),w)-2\ge 15|C|$.

Take
$p=[q/2]+1$
and consider the initial subsentence
$\al\sub\be$
with
$\cL(\al,w)=p+2$.
Then by the same reason as above, there are at least
$p$
stop signs behind
$a$
in
$\al$.

For the tail
$t_a(\al)$
of
$\al$
with initial letter
$a$,
we have
$|t_a(\al)|\le|t_a(\be)|$
and therefore
$$p\ge q/2\ge\cL(f_c(v),w)/2
  \ge\frac{|t_a(\be)|}{2|C|}-\frac{1}{2}
   \ge\frac{|t_a(\al)|}{3|C|}+2,$$
where we have used the estimate
$|t_a(\be)|\ge 15|C|$.
Let
$U\in T_c$
be the vertex (or covering element) corresponding
to the sentence
$\al$.
By Lemma~\ref{lem:critlevel}(2), there are at most three vertices
of the segment
$wU\sub wf_c(v)\sub T_c$
(including
$w$)
having the level
$<l$.
Because
$\cL(\al,w)=p+2\ge q/2+2\ge 15|C|/2+3\ge 9$,
we obtain
$\ell(U)\ge l+1$.

Assume that
$\cL(\be',w))\le\cL(\al,w)$.
Since
$\eta_c$
is radially isometric, we have
$\cL(\eta_c(\be),\eta_c(\be'))\ge\cL(\be,w)-\cL(\be',w)
  \ge\cL(\be,\al)$.
Because
$\cL(\be,\al)=q-p-2\ge q/2-3=\cL(f_c(v),w)/2-3$,
using estimate (2) above, we obtain
$$|\eta(v)\eta(v')|\ge\frac{|vv'|-\si}{2|C|}-3\ge\la|vv'|-\si'$$
with
$\la=1/2|C|$
and
$\si'=\si/2|C|+3\le 5$.

Now, we assume that
$\cL(\be',w))>\cL(\al,w)$.
Let
$U'\in T_c$
be the vertex corresponding to the initial subsentence
$\al'\sub\be'$
with
$\cL(\al',w)=\cL(\al,w)=p+2$.

There is no reason for the levels
$\ell(U)$, $\ell(U')$
to coincide. However, because
$\cL(\al,w)\ge 9$,
we have
$\ell(U')\ge l+1$
since there is at most three vertices of the segment
$wU'\sub wf_c(v')\sub T_c$
below the critical level.

Then, the condition of Proposition~\ref{pro:crlevelsentences}
is satisfied for the sentences
$\al=U$, $\al'=U'$
of the alphabet
$A$
because
$v\in f_c(v)\sub U$
and
$v'\in f_c(v')\sub U'$
by the definition of the map
$f_c$
and properties of the tree
$T_c$.
By this Proposition, the letters
$a\in(s_m,s_{m+1})\sub\be$,
$a'\in(s_{m'},s_{m'+1})\sub\be'$
of the critical level,
$\lv(a)=l=\lv(a')$,
are different,
$a\neq a'$,
and
$|m-m'|\le 2$.

Furthermore, the equality
$\cL(\al',w)=p+2$
means that there is at least
$p$
stop signs behind
$a'$
in
$\al'$.

For the tail
$t_{a'}(\be')$
of the sentence
$\be'$,
we have
$|t_{a'}(\be')|=\ell(f_c(v'))-l+1$.
Using (1) as above, we obtain
$$|t_{a'}(\be')|\le|C|(\cL(f_c(v),w)+1)
   =|C|(q+1)\le|C|(2p+1)\le 3|C|(p-2),$$
because
$p\ge 7$.
Thus,
$|t_{a'}(\al')|\le|t_{a'}(\be')|\le n(p-2)$
with
$n=3|C|$.

In this way, we have recovered the symmetry between
$\al$, $\al'$
needed to apply Proposition~\ref{pro:equaldiaries}:
there are at least
$p$
stop signs behind
$a$, $a'$
in
$\al$, $\al'$
respectively and
$\max\{|t_a|,|t_{a'}|\}\le n(p-2)$
for some
$p\ge 3$
and
$n=3|C|$.
Since the conditions (1)--(3) of that Proposition
are satisfied and
$a\neq a'$,
we conclude that
$\psi_c(\al)\neq\psi_c(\al')$
for the diary constant
$\ka\ge 15|C|+1$.

This yields the required estimate: because
$\psi_c$
is radially isometric, we have
$\cL(\psi_c(\be),\psi_c(\be'))\ge\cL(\be,\al)$
and as above we obtain
$$|\eta(v)\eta(v')|\ge\la|vv'|-\si'$$
with
$\la=1/2|C|$
and
$\si'\le 5$.
This completes the proof of Theorem~\ref{thm:compose}.
\end{proof}

\bigskip

\begin{tabbing}

Sergei Buyalo,\hskip11em\relax \= Viktor Schroeder,\\

St. Petersburg Dept. of Steklov \>
Institut f\"ur Mathematik, Universit\"at \\

Math. Institute RAS, Fontanka 27, \>
Z\"urich, Winterthurer Strasse 190, \\

191023 St. Petersburg, Russia\>  CH-8057 Z\"urich, Switzerland\\

{\tt sbuyalo@pdmi.ras.ru}\> {\tt vschroed@math.unizh.ch}\\

\end{tabbing}

\end{document}